\documentclass[a4paper,12pt,reqno]{amsart}
\usepackage{graphics}
\newtheorem{lem}{Lemma}
\newtheorem{con}{Conjecture}
\newtheorem{theo}{Theorem}

\newtheorem{prop}{Proposition}
\numberwithin{equation}{section}

\newcommand{\id}{\operatorname{id}}

\begin{document}

\title[Alternating sign matrices]{Linear relations of refined enumerations of alternating sign matrices}

\author[Ilse Fischer]{Ilse Fischer}

\thanks{Supported by the Austrian Science Foundation
    FWF, START grant Y463 and NFN grant S9607--N13.}

\begin{abstract}
In recent papers we have studied refined enumerations of alternating sign matrices with respect to a fixed set of top and bottom rows. The present paper is a first step towards extending these considerations to alternating sign matrices where in addition a number of left and right columns are fixed. The main result is a simple linear relation between the number of $n \times n$ alternating sign matrices where the top row as well as the left and the right column is fixed and the number of $n \times n$ alternating sign matrices where the two top rows and the bottom row is fixed. This may be seen as a first indication for the fact that the refined enumerations of alternating sign matrices with respect to a fixed set of top and bottom rows as well as left and right columns can possibly be reduced to the refined enumerations where only a number of top and bottom rows are fixed. For the latter numbers we provide a system of linear equations that conjecturally determines them uniquely. 
\end{abstract}

\maketitle

\section{Introduction}

The central objects of this article are {\it alternating sign matrices} which are defined as square matrices with entries 
in $\{1,-1,0\}$ such that in every row and column the following two conditions are fulfilled: the sum of entries is $1$ and the non-zero entries alternate in sign. 
For instance, 
$$\left( \begin{tabular}{rrrrr}
0 & 1 & 0 & 0 & 0 \\
0 & 0 & 1 & 0 & 0 \\
1 & -1 & 0 & 0 & 1 \\
0 & 1 & -1 & 1 & 0 \\
0 & 0 & 1 &  0 & 0
\end{tabular} \right)
$$
is a $5 \times 5$ alternating sign matrix. The fascination concerning these objects stems from the fact that they belong to the mysterious class of combinatorial objects with the property that their enumerations 
subject to a variety of different constraints lead to simple product formulas while at the same time proving these formulas is a non-trivial task. The first result in this respect was the enumeration of all alternating sign matrices of given size:  
Mills, Robbins and Rumsey \cite{mills1, mills2} conjectured and Zeilberger \cite{zeilberger} proved that there are 
$$
\prod_{j=0}^{n-1} \frac{(3j+1)!}{(n+j)!}=:A_n
$$
$n \times n$ alternating sign matrices. We refer to Bressoud's beautiful book \cite{bressoud} for an introduction into this field. Later it turned out that also the enumeration of many {symmetry classes of alternating sign matrices} of fixed size leads to enumeration formulas of this type. 

\medskip

\subsection{Some known results on linear relations between refined enumerations of alternating sign 
matrices} This article is concerned with refined enumerations of alternating sign matrices. The first result in this respect 
was also already discovered by Mills, Robbins and Rumsey \cite{mills2} and later proved by Zeilberger \cite{zeilberger2}: it is not hard to see that an alternating sign matrix has a unique $1$ in its top row. Surprisingly, the number of $n \times n$ alternating sign matrices with a $1$ in column $k$ of the top row can again be expressed by a product formula, namely 
$$
\binom{n+k-2}{k-1} \frac{(2n-k-1)!}{(n-k)!} \prod_{j=0}^{n-2} \frac{(3j+1)!}{(n+j)!} =:A_{n,k}.
$$
In \cite{newproof}, I gave an alternative proof of this formula. It
is based on the (non-trivial) observation that the numbers $A_{n,k}$ are a solution of the following system 
of linear equations  
\begin{equation}
\label{binomrelation1}
\sum_{i=1}^{n} A_{n,i} (-1)^{i+1} \binom{2n-k-1}{n-k-i+1} = A_{n,k}, \qquad 1 \le k \le n.
\end{equation}
Namely, it turned out that the solution space of this system is one-dimensional, and, consequently, determines the numbers $A_{n,k}$ uniquely together with the obvious identity $A_{n,1} = A_{n-1} = 
\sum\limits_{k=1}^{n-1} A_{n-1,k}$ inductively with $n$.  Thus, in order to 
complete the proof of the so-called refined alternating sign matrix theorem, it was essentially sufficient to show that 
the conjectural formula for $A_{n,k}$ provides a solution for the system of linear equations. This puts everything down to certain 
hypergeometric identities. (Usually it is easy to come up with a 
conjectural formula in this type of enumeration problems. The art really lies in proving the formula.)
In the present article we deal with linear relations between different types of more refined enumerations of alternating sign matrices.

\medskip

There exist even simpler linear relations that were previously discovered: in our next example
two different types of doubly refined enumerations of alternating sign matrices are involved. Stroganov \cite{stroganov} considered the numbers $\overline{\underline{A}}_{n,i,j}$ of $n \times n$ 
alternating sign matrices with a $1$ in position $(1,i)$ and in position $(n,j)$ and the numbers $|\overline{A}_{n,i,j}$ of 
$n \times n$ alternating sign matrices with $1$s in position $(1,i)$ and in position $(j,1)$. It is certainly too much 
to expect also product formulas for these numbers (Stroganov implicitly gave an explicit formula for these numbers in terms of a sum of products, see, for instance, \cite{fischertrapezoids}), but Stroganov showed that they are 
linearly related as follows
\begin{equation}
\label{stroganov}
\overline{\underline{A}}_{n,i,j} = |\overline{A}_{n,i+1,j} + |\overline{A}_{n,i,j+1} - |\overline{A}_{n,i+1,j+1} =
(E_i + E_j - E_i E_j) |\overline{A}_{n,i,j}
\end{equation}
if $i,j \in \{1,2,\ldots,n\}$ and $i+j \ge 3$, and where $E_i$ denotes the shift operator, i.e. $E_i a(i) = a(i+1)$. 
(If $i=j=1$ then $A_{n,i,j} = [n=1]$, where $[\text{statement}]=1$ if the statement is true and 
$[\text{statement}]=0$ otherwise.) As it is the standard situation for phenomenons related to alternating sign matrices also this simple fact does not seem to have a simple (bijective) explanation so far.

\medskip

Another doubly refined enumeration of alternating sign matrices was considered in \cite{fischerromik}: if the first two rows of an $n \times n$ alternating sign matrix are deleted then we obtain a matrix where in exactly two columns, say columns $i$ and $j$, the first non-zero entry is a $-1$. 
Let $\overline{\overline{A_{n,i,j}}}$ ($i < j$) denote the number of these matrices. All alternating sign 
matrices that lead to such a ``partial'' $(n-2) \times n$ alternating sign matrix can be obtained by choosing an integer $k$ with $i \le k \le j$ and adding $e_k:=([l=k])_{1 \le l \le n}$ as first row and $e_i+e_j-e_k$ as 
second row. In this sense, the numbers $\overline{\overline{A_{n,i,j}}}$ enumerate alternating sign matrices with respect to the two top rows. Karklinsky and Romik \cite{2topbottom} were able to derive the following linear relation 
between the numbers $\overline{\overline{A_{n,i,j}}}$ and  $\overline{\underline{A}}_{n,i,j}$ 
\begin{equation}
\label{kark}
\overline{\underline{A}}_{n,i+1,j+1} - \overline{\underline{A}}_{n,i,j} = 
\overline{\overline{A_{n,i+1,n+1-j}}} + \overline{\overline{A_{n,i,n-j}}} - \overline{\overline{A_{n,i,n+1-j}}} -
\overline{\overline{A_{n,n-j,i}}} - \overline{\overline{A_{n,n+1-j,i+1}}} + \overline{\overline{A_{n,n-j,i+1}}}
\end{equation}
(we set $\overline{\underline{A}}_{n,i,j}=0$ and $\overline{\overline{A_{n,i,j}}}
=0$ if $i$ and $j$ are outside of the respective range of their definition) and used this and Stroganov's formula for $\overline{\underline{A}}_{n,i,j}$ to deduce an explicit formula 
for $\overline{\overline{A_{n,i,j}}}$. Interestingly, in \cite{fischertrapezoids}, another linear relation between these numbers was derived, namely 
\begin{equation}
\label{binomrelation2}
\overline{\overline{A_{n,i,j}}} = \sum_{k=j}^{n} (-1)^{n+k} \binom{2n-2-j}{k-j} \overline{\underline{A}}_{n,i,k}
\end{equation}
for $1 \le i < j \le n$, which was used to deduce another formula (of the same type) for $\overline{\overline{A_{n,i,j}}}$.

\medskip

With so many linear relations between doubly refined enumeration numbers around it is natural to ask for linear 
relations between triply refined enumeration numbers which possibly generalize the former. We introduce the three triply refined enumerations of alternating sign matrices that we will consider. To this end, we need a little background from \cite{fischertrapezoids}, where we have studied refined enumerations of 
alternating sign matrices with respect to the $d$ top rows and the $c$ bottom rows.

\medskip

\subsection{Refined enumerations of alternating sign matrices with respect to the $d$ top rows and the $c$ bottom rows} In order to explain roughly what we have accomplished there, we fix an $n \times n$ alternating sign matrix $A$ and two non--negative integers $c, d$ with $c+d \le n$. Let $T$ denote the $d \times n$ matrix, which consists of the $d$ top rows of $A$, let $B$ denote the $c \times n$ matrix which consists the $c$ bottom rows of $A$ and $M$ denote the $(n-c-d) \times n$ matrix which we obtain from $A$ after we have deleted $T$ and $B$.  

\medskip

Clearly, in the three matrices $T$, $M$ and $B$, the non-zero entries alternate in each row and column and the row sums are $1$. However, there are precisely $d$ columns in $T$, say columns $i_1, i_2, \ldots, i_d$ with $1 \le i_1 < i_2 < \ldots < i_d \le n$, whose column sums are $1$. (If we choose $d=2$ in the example above then 
$(i_1,i_2)=(2,3)$.) The other columns of $T$ sum up to $0$ and their first non-zero entry (which may or may not exist) is a $1$. These properties characterize the matrices $T$ which may appear in this way and, in fact, they correspond to monotone triangles with $d$ rows and bottom row $(i_1,i_2,\ldots,i_d)$. (For the definition of monotone triangles see Section \ref{prel}.) Let $\alpha(d;i_1,\ldots,i_d)$ 
denote the number of these objects. (Note that this number does not depend on $n$ as the first $i_1-1$ columns as well as the last $n-i_d$ columns of $T$ consist entirely of zeros.)

\medskip

Likewise, there exist $c$ columns in $B$, say columns $s_1, s_2, \ldots, s_c$ with $1 \le s_1 < s_2 < \ldots < s_c \le n$, whose column sums are $1$. (If we choose $c=2$ in the example then $(s_1,s_2)=(2,4)$.) The other columns of $B$ sum up to $0$ and their first non-zero entry (if it exists) is a $-1$. Again, these properties characterize the matrices $B$  and, by reflection along a horizontal axis, the number of the possible matrices $B$ is equal to $\alpha(c;s_1,\ldots,s_c)$.

\medskip

The 
matrices $M$ are characterized by the following properties (in addition to the facts that the entries alternate in rows and columns and that all row sums are $1$): the column sum is $1$ precisely for the columns whose index is neither in $\{i_1,\ldots,i_d\}$ nor 
in $\{s_1,\ldots,s_c\}$, the column sum is $-1$ precisely for the columns whose index is both in $\{i_1,\ldots,i_d\}$ and in $\{s_1,\ldots,s_c\}$. For all other columns, the sum is zero and the first non-zero entry (if it exists) is a $-1$ if and only if the column is in $\{i_1,\ldots,i_d\}$. In the following, we let $A(n;s_1,\ldots,s_c;i_1,\ldots,i_d)$ denote the number of these matrices $M$. Observe also that conversely, in order to obtain all $n \times n$ alternating sign matrix, we choose any two strictly increasing sequences $(i_1,\ldots,i_d)$ and $(s_1,\ldots,s_c)$ of integers in 
$\{1,2,\ldots,n\}$, and combine any triple $(T,M,B)$ of matrices with the properties given above. Phrased differently, 
$$
A_n = \sum_{1 \le s_1 < s_2 < \ldots < s_c \le n \atop 1 \le i_1 < i_2 < \ldots < i_d \le n} \alpha(d;i_1,\ldots,i_d) A(n;s_1,\dots, s_c;i_1,\ldots,i_d) \alpha(c;s_1,\ldots,s_c).
$$

\medskip

To make a long story short, the numbers $A(n;s_1,\ldots,s_c;i_1,\ldots,i_d)$  obviously enumerate $n \times n$ alternating sign matrices with respect to a fixed set of $d$ top rows and a fixed set of 
$c$ bottom rows. Some special cases have already appeared above: $A_n = A(n;-;-)$, $A_{n,i} = A(n;i;-) = A(n;-;i)$, $\overline{\underline{A}}_{n,i,j} = A(n;i;j) = A(n;j;i)$ and $\overline{\overline{A}}_{n,i,j} = A(n;i,j;-) = A(n;-;i,j)$. In Theorem~1 of \cite{fischertrapezoids} (see also Section~\ref{prel}) we have identified 
these refined enumeration numbers $A(n;s_1,\ldots,s_c;i_1,\ldots,i_d)$ as the coefficients of a specialization of $\alpha(n;k_1,\ldots,k_n)$ with respect to a certain polynomial basis in the parameters $k_1,\ldots,k_n$. (It turns out that $\alpha(n;k_1,\ldots,k_n)$ is a polynomial in $(k_1,\ldots,k_n)$. We refer to Section~\ref{prel} for details on $\alpha(n;k_1,\ldots,k_n)$.)

\medskip

\subsection{Triply refined enumerations of alternating sign matrices -- the main theorem} Among the  numbers $A(n;s_1,\ldots,s_c;i_1,\ldots,i_d)$, there are two essential 
different triply refined enumerations of alternating sign matrices: on the one hand, this is the quantity $A(n;-;i_1,i_2,i_3)=:\overline{\overline{\overline{A_{n,i_1,i_2,i_3}}}}$, which amounts to enumerate 
$n \times n$ alternatig sign matrices with respect to the three top rows. On the other hand, we have $A(n;s_1;i_1,i_2)=:\overline{\overline{\underline{A_{n,s_1,i_1,i_2}}}}$, which amounts to enumerate 
$n \times n$ alternating sign matrices with respect to the two top rows and with respect to the bottom row.

\medskip

In \cite{fischertrapezoids}, it also turned out 
that all the refined enumeration numbers $A(n;s_1,\ldots,s_c;i_1,\ldots,i_d)$ where $c$ and $d$ sum up to the same integer are linearly related. To be more accurate, we have 
\begin{multline}
\label{topbottom1}
A(n;s_1,\ldots,s_c;i_1,\ldots,i_d) = \sum_{i_{d+1}=s_c}^{n} \sum_{i_{d+2}=s_{c-1}}^{n} \cdots 
\sum_{i_{d+t}=s_{c-t+1}}^{n} A(n;s_1,\ldots,s_{c-t};i_1,\ldots,i_{d+t}) \\
\times (-1)^{i_{d+1}+\ldots+i_{d+t}+t n} \binom{2n-c-d-s_c}{i_{d+1}-s_c} \binom{2n-c-d-s_{c-1}}{i_{d+2}-s_{c-1}} 
\cdots \binom{2n-c-d-s_{c-t+1}}{i_{d+1} - s_{c-t+1}}
\end{multline}
for $1 \le t \le c$ (this is identity (5.1) from \cite{fischertrapezoids}) as well as
\begin{multline}
\label{topbottom2}
A(n;s_1,\ldots,s_c;i_1,\ldots,i_d) = \sum_{s_{c+1}=i_d}^{n} \sum_{s_{c+2}=i_{d-1}}^{n} \cdots 
\sum_{s_{c+t}=i_{d-t+1}}^{n} A(n;s_1,\ldots,s_{c+t};i_1,\dots,i_{d-t}) \\
\times (-1)^{s_{c+1}+\ldots+s_{c+t}+t n} \binom{2n-c-d-i_{d}}{s_{c+1}-i_d} \binom{2n-c-d-i_{d-1}}{s_{c+2}-i_{d-1}} 
\cdots \binom{2n-c-d-i_{d-t+1}}{s_{c+t}-i_{d-t+1}}
\end{multline}
for $1 \le t \le d$. These identities in fact generalize \eqref{binomrelation1} and \eqref{binomrelation2}. However, note that they involve the quantity $A(n;s_1,\ldots,s_c;i_1,\ldots,i_d)$ 
also for sequences $(s_1,\ldots,s_c) \in \{1,\ldots,n\}^c$ and $(i_1,\ldots,i_d) \in \{1,2,\ldots,n\}^d$ that 
are not strictly increasing, which have (so far) no combinatorial meaning. In these
``non-combinatorial'' cases, the definition of the numbers is via their interpretation
as certain coefficients of a specialization of $\alpha(n;k_1,\ldots,k_n)$, see 
\cite[Theorem~1]{fischertrapezoids} or Section~\ref{prel}. (Note that from this point of 
view we have $A_{n,i,j} = A(n;i,j;-) \not= 0$ for $i \ge j$ if $i, j \in \{1,2,\ldots,n\}$. This is 
in contrary to the situation in \eqref{kark} and the viewpoint we take for the rest of the paper.)

\medskip

As a consequence of these identities, we see that the two triply refined enumerations $\overline{\overline{\overline{A_{n,i_1,i_2,i_3}}}}$ and $\overline{\overline{\underline{A_{n,s_1,i_1,i_2}}}}$ are linearly related as follows. Equation~\eqref{topbottom1} implies
\begin{equation}
\label{321}
\overline{\overline{\underline{A_{n,s_1,i_1,i_2}}}} = \sum_{i_3=s_1}^{n} \overline{\overline{\overline{A_{n,i_1,i_2,i_3}}}} (-1)^{i_3+n} \binom{2n-3-s_1}{i_3-s_1}
\end{equation}
and Equation~\eqref{topbottom2} implies
\begin{equation*}
\overline{\overline{\overline{A_{n,i_1,i_2,i_3}}}} = \sum_{s_1=i_3}^{n} \overline{\overline{\underline{A_{n,s_1,i_1,i_2}}}} (-1)^{s_1+n} \binom{2n-3-i_3}{s_1-i_3}.
\end{equation*}
(Unfortunately, for $s_1 \le i_2$, the right-hand side of the first identity involves numbers $\overline{\overline{\overline{A_{n,i_1,i_2,i_3}}}}$ that have no combinatorial interpretation. On the other hand, the right-hand side of the second identity involves only numbers that have a combinatorial interpretation even if $i_3 \le i_2$, in which case the left-hand side has no combinatorial interpretation.)

\medskip

However, there is of course a third (and maybe most obvious) possibility for a triply refined enumeration of alternating sign matrices: the enumeration of $n \times n$ alternating sign matrices with respect to 
the top row, the leftmost column and the rightmost column. In the following, $|\overline{A}|_{s,i,t}$ denotes the number of 
$n \times n$ alternating sign matrices with $1$s in positions $(s,1)$, $(1,i)$ and $(t,n)$. The main result of the present paper is the following simple linear relation between the numbers 
$\overline{\overline{\underline{A_{n,s,i,t}}}}$ and the numbers $|\overline{A}|_{n,s,i,t}$. The most fascinating fact concerning this paper is probably that this theorem once more constitutes a result on alternating sign matrices that can be stated easily, yet it seems to require a rather complicated proof. I would be eager to see a simpler proof, preferable more combinatorial and less algebraic.

\medskip

\begin{theo}
\label{main}
For $s,t \in \{2,3,\ldots,n\}$ and $i \in \{1,2,\ldots,n\}$, we have 
\begin{multline}
\label{triplytheo}
\overline{\overline{\underline{A_{n,s,i,t}}}} = |\overline{A}|_{n,s+1,i,t}
+ |\overline{A}|_{n,s+1,i,t+1} + |\overline{A}|_{n,s,i+1,t} 
- |\overline{A}|_{n,s+1,i+1,t} - |\overline{A}|_{n,s,i,t+1} 
- |\overline{A}|_{n,s+1,i-1,t+1} \\
= (E_s + E_s E_t + E_i - E_s E_i - E_t - E_s E^{-1}_i E_t) |\overline{A}|_{n,s,i,t}
\end{multline}
where we set $|\overline{A}|_{n,n+1,i,t} = |\overline{A}|_{n,s,i,n+1} = |\overline{A}|_{n,s,0,t} =
|\overline{A}|_{n,s,n+1,t}=0$.
\end{theo}

\medskip

Stroganov's identity \eqref{stroganov} is now a corollary of this, which we obtain if we specialize $t=n$.  The proof of Theorem~\ref{main} is presented in Section~\ref{triply}, after we have summarized necessary preliminaries in Section~\ref{prel} and derived (as a warm-up) another linear relation between two types of doubly refined enumerations of alternating sign matrices in Section~\ref{warmup}. 

\medskip

The approach we chose to prove Theorem~\ref{main} heavily uses various properties of the polynomial $\alpha(n;k_1,\ldots,k_n)$. We believe that the six-vertex model approach should also be applicable to derive it for the combinatorial meaningful cases, i.e. if $i < t$.  (See, for instance, \cite{bressoud} to learn about the six--vertex model.) In this case, the symmetry of the partition function in the 
row and column ``spectral'' parameters should play a fundamental role, see also the derivation of \eqref{kark} in \cite{2topbottom}. However, we do not see how the non-combinatorial cases (i.e. if $i \ge t$) can be dealt with using this approach. 

\medskip

\subsection{A system of linear equations for the refined enumeration of alternating sign matrices with respect to the $d$ top rows -- the main conjecture} Since we now know that the three types of triply refined enumerations of alternating sign matrices are linearly related it would be desirable to find a possibility to compute either of these numbers. 

\medskip

Indeed, the relation in Theorem~\ref{main}
can also be used to express the numbers $|\overline{A}|_{n,s,i,t}$ in terms of the numbers $\overline{\overline{\underline{A_{n,s,i,t}}}}$ (which finally implies that any triply refined enumeration can be expressed in terms of any other): in order to see that observe that for $n \ge 2$, we have the following behavior of $|\overline{A}|_{n,s,i,t}$ on the boundary of the combinatorial admissible domain, i.e. for 
$s=1, n$ and $t=1, n$:
$$
|\overline{A}|_{n,1,i,t} = [t \not= 1] [i=1] A_{n-1,t-1}, \qquad 
|\overline{A}|_{n,s,i,1} = [s \not= 1] [i=n] A_{n-1,s-1},
$$
$$
|\overline{A}|_{n,n,i,t} = [t \not= n] [i \not= 1] \overline{A}|_{n-1,i-1,t} \quad  \text{and} \quad
|\overline{A}|_{n,s,i,n} = [s \not= n] [i \not= n] |\overline{A}_{n-1,s,i}.
$$
An explicit formula for $A_{n,k}$ is given above and formulas for $|\overline{A}_{n,s,i}$ as well as for $\overline{A}|_{n,i,t}$ are given in Section~\ref{warmup}.
If we rewrite the identity in Theorem~\ref{main} as follows
$$
 |\overline{A}|_{n,s,i+1,t}  = \overline{\overline{\underline{A_{n,s,i,t}}}}  -  |\overline{A}|_{n,s+1,i,t}
- |\overline{A}|_{n,s+1,i,t+1}  
+ |\overline{A}|_{n,s+1,i+1,t} + |\overline{A}|_{n,s,i,t+1} 
+ |\overline{A}|_{n,s+1,i-1,t+1}
$$
then the number $|\overline{A}|_{n,s,i+1,t}$ is expressed in terms of $\overline{\overline{\underline{A_{n,s,i,t}}}}$ and  several $|\overline{A}|_{n,s',i',t'}$'s with either 
$s' =s+1 > s$ or $t' =t+1 > t$. This also shows why the non-combinatorial 
cases of the relation in Theorem~\ref{main} are of importance rather than a secondary result: when expressing $|\overline{A}|_{n,s,i,t}$ in terms of the numbers $\overline{\overline{\underline{A_{n,s,i,t}}}}$ then, in general, also instances of the latter numbers where $i \ge t$ are involved.

\medskip

We present a system of linear equations that conjecturally determines the numbers $\overline{\overline{\overline{A_{n,i_1,i_2,i_3}}}}=:A(n;-;i_1,i_2,i_3)$ uniquely. In fact, this seems to extend to the numbers $A(n;-;i_1,\ldots,i_d)$ for $d \ge 4$ and also to $d=1,2$. The relation in \eqref{topbottom1} can then be used to compute $A(n;s_1,\ldots,s_c;i_1,\ldots,i_d)$ also for the cases where $c > 0$.

\begin{con} 
\label{sufficiency}
Let $n, d$ be integers with $n \ge d \ge 1$ and assume that we have already computed 
the numbers $A(n-1;-;i_1,\ldots,i_{d-1})$ for $1 \le i_1 < i_2 < \ldots < i_{d-1} \le n-1$. Then the numbers $A(n;-;i_1,\ldots,i_d)$, $1 \le i_1, i_2, \ldots,i_d \le n$, are uniquely determined by the following system of linear equations.
\begin{enumerate}
\item For all $i_1,\ldots,i_d \in \{1,2,\ldots,n\}$, 
\begin{multline*}
A(n;-;i_1,\ldots,i_d) = \sum_{j_1=i_1}^{n} \sum_{j_2=i_2}^{n} \dots \sum_{j_d=i_d}^{n} (-1)^{d n + j_1 + \ldots + j_d} 
A(n;-;j_d,\ldots,j_1) \\ \times \binom{2n-i_1-d}{j_1-i_1} \binom{2n-i_2-d}{j_2-i_2} \cdots \binom{2n-i_d-d}{j_d-i_d}.
\end{multline*}
\item For all $l \in \{1,2,\ldots,d-1\}$ and $i_1,\ldots,i_d \in \{1,2,\ldots,n\}$,
\begin{multline*}
A(n;-;i_1,\ldots,i_d) - A(n;-;i_1,\ldots,i_l,i_{l+1}+1,i_{l+2},\ldots,i_d) 
+ A(n;-;i_1,\ldots,i_l+1,i_{l+1}+1,\ldots,i_d) \\
+ A(n;-;i_1,\ldots,i_{l-1},i_{l+1},i_l,i_{l+2},\ldots,i_d) 
-  A(n;-;i_1,\ldots,i_{l-1},i_{l+1},i_l+1,i_{l+2},\ldots,i_d) \\
+  A(n;-;i_1,\ldots,i_{l-1},i_{l+1}+1,i_l+1,i_{l+2},\ldots,i_d) = 0
\end{multline*}
\item For $1 \le i_1 < i_2 < \ldots < i_d \le n$, 
$$
A(n;-;i_1,\ldots,i_d) = A(n;-;n+1-i_d,n+1-i_{d-1},\ldots,n+1-i_1).
$$
\item For $1 \le i_1 < i_2 < \ldots < i_{d-1} \le n-1$, 
$$
A(n;-;i_1,\ldots,i_{d-1},n) = A(n-1;-;i_1,\ldots,i_{d-1}).
$$
\end{enumerate}
\end{con}

\medskip

The fact that the numbers $A(n;-;i_1,\ldots,i_d)$ are a solution to this system of linear equations is (more or less) already proven in previous papers. To be more accurate:
\begin{enumerate}
\item The first set of equations was conjectured in \cite[Conjecture~5]{fischerromik} and proven in \cite[Equation (6.1)]{fischertrapezoids}.
\item The second set of equations follows from  the translation of the operator identity
\eqref{swap} for $\alpha(n;k_1,\ldots,k_n)$ to its coefficients, see \eqref{expansion}. The crucial fact is $\delta_{k_{n-d+l}} \binom{k_{n-d+l}-n+d-2+i_{l}}{i_{l}-1} = \binom{k_{n-d+l}-n+d-2+i_l-1}{i_l-2}$, where $\delta_x = \id - E^{-1}_x$. 
\item The third and the fourth set of equations follows from the combinatorial interpretation of the numbers 
$A(n;-;i_1,\ldots,i_d)$.
\end{enumerate}
Thus, in order to prove Conjecture~\ref{sufficiency}, it remains to solve the basic linear algebra exercise that the displayed systems of linear equations have at most 
one solution.  Finally, it should be noted that there is also another possibility to compute the numbers $A(n;-;i_1,\ldots,i_d)$: let $(k_1,\ldots,k_{n-d})$ be the 
strictly increasing sequence of integers with $\{i_1,\ldots,i_d\} \cup \{ k_1, \ldots, k_{n-d}\} = \{1,2, \ldots, n\}$. Then, by the definitions of $A(n;s_1,\ldots,s_c;i_1,\ldots,i_d)$ and $\alpha(n;k_1,\ldots,k_n)$, we have 
$A(n;-;i_1,\ldots,i_d) = \alpha(n-d;k_1,\ldots,k_{n-d})$. An explicit formula for $\alpha(n;k_1,\ldots,k_n)$ is provided in \eqref{operator} below.

\section{Preliminaries: some important facts about $\alpha(n;k_1,\dots,k_n)$}
\label{prel}

In this section we summarize several notions and results from previous papers. 

\medskip

\subsection{Monotone triangles} 
A {\it monotone triangle} is an integer array of the following triangular shape 
\begin{center}
\begin{tabular}{ccccccccccccccccc}
  &   &   &   &   &   &   &   & $a_{n,n}$ &   &   &   &   &   &   &   & \\
  &   &   &   &   &   &   & $a_{n-1,n-1}$ &   & $a_{n-1,n}$ &   &   &   &   &   &   & \\
  &   &   &   &   &   & $\dots$ &   & $\dots$ &   & $\dots$ &   &   &   &   &   & \\
  &   &   &   &   & $a_{3,3}$ &   & $\dots$ &   & $\dots$ &   & $a_{3,n}$ &   &   &   &   & \\
  &   &   &   & $a_{2,2}$ &   & $a_{2,3}$ &  &   $\dots$ &   & $\dots$   &  & $a_{2,n}$  &   &   &   & \\
  &   &   & $a_{1,1}$ &   & $a_{1,2}$ &   & $a_{1,3}$ &   & $\dots$ &   & $\dots$ &   & $a_{1,n}$ &   &   &
\end{tabular},
\end{center}
which is monotone increasing in northeast and in southeast
direction and strictly increasing along rows, that is $a_{i,j} \le
a_{i+1,j+1}$ for $1 \le i \le j < n$, $a_{i,j} \le a_{i-1,j}$ for
$1 < i \le j \le n$ and $a_{i,j} < a_{i,j+1}$ for $1 \le i \le j
\le n-1$. Monotone triangles with bottom row $(1,2,\ldots,n)$ are said to be complete and correspond to 
$n \times n$ alternating sign matrices: for a given complete monotone triangle $A=(a_{i,j})_{1 \le i \le j \le
n}$ with $n$ rows,  the corresponding $n \times n$
alternating sign matrix $M=(m_{i,j})_{1 \le i, j \le n}$ can be obtained as follows: we have $m_{i,j}=1$
if and only if
$$j \in
\{a_{n+1-i,n+1-i},a_{n+1-i,n+2-i},\ldots,a_{n+1-i,n}\} \setminus
\{a_{n+2-i,n+2-i},a_{n+2-i,n+3-i},\ldots,a_{n+2-i,n}\},
$$
$m_{i,j} = -1$ if and only if
$$j \in
\{a_{n+2-i,n+2-i},a_{n+2-i,n+3-i},\ldots,a_{n+2-i,n}\} \setminus
\{a_{n+1-i,n+1-i},a_{n+1-i,n+2-i},\ldots,a_{n+1-i,n}\}
$$
and $m_{i,j}=0$ else.

\medskip

\subsection{The operator formula and other properties of $\alpha(n;k_1,\ldots,k_n)$} 
Let $\alpha(n;k_1,\ldots,k_n)$ denote the number of monotone triangles with bottom row 
$(k_1,\ldots,k_n)$, i.e. the integer arrays $(a_{i,j})_{1 \le i \le j \le n}$ with the properties 
given above and with
$a_{1,i}=k_i$. (Note that this is compatible with the definition of $\alpha$ given above.) This quantity and its numerous properties lie in the ``heart'' of our proofs. 
To begin with, $\alpha$ can obviously be computed recursively as follows.
\begin{multline}
\label{recursion}
\alpha(n;k_1,\ldots,k_n) = \sum_{(l_1,\ldots,l_{n-1}) \in \mathbb{Z}^{n-1} \atop 
k_1 \le l_1 \le k_2 \le l_2 \le k_3 \le \ldots \le k_{n-1} \le l_{n-1} \le k_n, l_i \not= l_{i+1}}
\alpha(n-1;l_1,\ldots,l_{n-1}) \\
=: \sum_{(l_1,\ldots,l_{n-1})}^{(k_1,\ldots,k_n)} \alpha(n-1;l_1,\ldots,l_{n-1})
\end{multline}
The following formula for this quantity was derived in \cite{fischertriangle,fischergtriangle} 
\begin{equation}
\label{operator}
\alpha(n;k_1,\ldots,k_n) = \prod_{1 \le p < q \le n} (\id - E_{k_p} + E_{k_p} E_{k_q}) \prod_{1 \le i < j \le n} 
\frac{k_j - k_i}{j-i} , 
\end{equation}
where $E_x$ denotes the shift operator, defined as $E_x p(x) = p(x+1)$, and ``$\id$'' is the identity operator. (Note that this expression for $\alpha$ extends its definition to $(k_1,\ldots,k_n) \in \mathbb{Z}^n$ that are not necessarily strictly increasing.)
The formula implies that 
\begin{equation*}
(\id  + S_{k_i,k_{i+1}}) (\id - E_{k_{i+1}} + E_{k_i} E_{k_{i+1}}) \alpha(n;k_1,\ldots,k_n) = 0,
\end{equation*}
where $S_{x,y}$ denotes the operator that swaps the variables $x$ and $y$, see \cite{fischertriangle}. This is equivalent to
\begin{equation}
\label{swap}
(\id  + S_{k_i,k_{i+1}}) (\id - \delta_{k_{i+1}} + \delta_{k_i} \delta_{k_{i+1}}) \alpha(n;k_1,\ldots,k_n) = 0,
\end{equation}
where $\delta_x = \id - E_x^{-1}$ as $E^{-1}_{k_i} E^{-1}_{k_{i+1}}$ commutes with
$S_{k_i,k_{i+1}}$  and 
$$
E^{-1}_{k_i} E^{-1}_{k_{i+1}} (\id - E_{k_{i+1}} + E_{k_i} E_{k_{i+1}}) = 
(\id - \delta_{k_{i+1}} + \delta_{k_i} \delta_{k_{i+1}}).
$$

\medskip

Another identity for $\alpha(n;k_1,\ldots,k_n)$, which is taken from Lemma~5 in \cite{newproof}, is also quite useful.
\begin{equation}
\label{cycle}
\alpha(n;k_1,\ldots,k_n) = (-1)^{n-1} \alpha(n;k_2,\ldots,k_n,k_1-n)
\end{equation}

\medskip

\subsection{A useful lemma on $\alpha(n;k_1,\ldots,k_n)$} 
The operator formula implies that $\alpha(n;k_1,\ldots,k_n)$ is a polynomial in $k_i$ of degree no greater than 
$n-1$. Thus, for every fixed $r \in \{1,2,\ldots,n\}$ and every integer $q$, it has to have a unique expansion in terms of the polynomial basis $\left( \binom{k_r + n - 1 - p - q}{n-1} \right)_{0 \le p \le n-1}$. We determine this expansion 
in the next lemma. (To see this, set $z=k_r-q$ and $k_r=q$ in the statement of the lemma.)
For this purpose, we need the following fact about $\alpha$: let $e_p(X_1,\ldots,X_n)$ denote the $p$--th elementary symmetric function. In \cite[Remark~1]{newproof}, it was shown that
\begin{equation}
\label{123}
e_p(\Delta_{k_1},\ldots,\Delta_{k_n}) \alpha(n;k_1,\ldots,k_n) = 0
\end{equation}
for $p \ge 1$, where $\Delta_x := E_x - \id$ denotes the difference operator. 

\medskip

\begin{lem}
\label{step3}
Let $1 \le r \le n$. Then
\begin{multline*}
\alpha(n;k_1,\ldots,k_{r-1},k_{r}+z,k_{r+1},\ldots,k_n)  \\
= \sum_{p=0}^{\infty} (-1)^p \binom{n-1-p+z}{n-1}
e_p(E_{k_1},\ldots,\widehat{E_{k_r}},\ldots,E_{k_n}) \, \alpha(n;k_1,\ldots,k_n)
\end{multline*}
where $\widehat{E_{k_r}}$ indicates that $E_{k_r}$ is missing in the argument.
\end{lem}

{\it Proof.} 
First observe that the left-hand side of the statement  is equal to
\begin{equation}
\label{z3}
E_{k_{r}}^{z} \alpha(n;k_{1},\ldots,k_{n})
= (\id + \Delta_{k_{r}})^{z} \alpha(n;k_{1},\ldots,k_{n})
= \sum_{i=0}^{z} \binom{z}{i} \Delta_{k_{r}}^{i}
\alpha(n;k_{1},\ldots,k_{n}),\end{equation}
where $\Delta_x: = E_x - \id$ denotes the difference operator.
We set
$$
\alpha_{i,r}(n;k_{1},\ldots,k_{n}) :=
e_{i}(\Delta_{k_{1}},\ldots,\widehat{\Delta_{k_{r}}},\ldots,\Delta_{k_{n}})
\, \alpha(n;k_{1},\ldots,k_{n}).
$$
By
\begin{equation}
\label{elementar-sym3}
X_{r} e_{i-1}(X_{1},\ldots,\widehat{X_{r}},\ldots,X_{n}) =
e_{i}(X_{1},\ldots,X_{n}) - e_{i}(X_{1},\ldots,\widehat{X_{r}},\ldots,X_{n})
\end{equation}
and \eqref{123},  we have
$$
\Delta_{k_{r}} \alpha_{i-1,r}(n;k_{1},\ldots,k_{n}) =
- \alpha_{i,r}(n;k_{1},\ldots,k_{n})
$$
for $i \ge 1$. This implies
\begin{equation}
\label{combine3}
\Delta_{k_{r}}^{i} \alpha_{0,r}(n;k_{1},\ldots,k_{n}) = (-1)^{i} \alpha_{i,r}(n;k_{1},\ldots,k_{n})
\end{equation}
for $i \ge 0$. As $\alpha_{0,r}(n;k_1,\ldots,k_n)=\alpha(n;k_1,\ldots,k_n)$, this implies
$$
E_{k_{r}}^{z} \alpha(n;k_{1},\ldots,k_{n})
= \sum_{i=0}^{z} (-1)^{i} \binom{z}{i}
e_{i}(\Delta_{k_{1}},\ldots,\widehat{\Delta_{k_{r}}},\ldots,\Delta_{k_{n}})
\alpha(n;k_{1},\ldots,k_{n}).
$$
Next we use the fact that
\begin{multline*}
e_{i}(\Delta_{k_{1}},\ldots,\widehat{\Delta_{k_{r}}},\ldots,\Delta_{k_{n}}) =
e_{i}(E_{k_{1}} - \id,\ldots,\widehat{E_{k_{r}} - \id},\ldots,E_{k_{n}} - \id) \\
= \sum_{p=0}^{\infty} \binom{n-1-p}{i-p} (-1)^{i+p}
e_{p}(E_{k_{1}},\ldots,\widehat{E_{k_{r}}},\ldots,E_{k_{n}})
\end{multline*}
in order to see that
\begin{multline*}
E_{k_{r}}^{z} \alpha(n,m;k_{1},\ldots,k_{n})  \\
= \sum_{i=0}^{z} \sum_{p=0}^{\infty} (-1)^{p}
\binom{n-1-p}{i-p} \binom{z}{i}
e_{p}(E_{k_{1}},\ldots,\widehat{E_{k_{r}}},\ldots,E_{k_{n}})
\alpha(n;k_{1},\ldots,k_{n}).
\end{multline*}
The assertion now follows from the Chu-Vandermonde summation. \qed

\medskip

\subsection{The numbers $A(n;s_1,\ldots,s_c;i_1,\ldots,i_d)$ appear as the coefficients of a specialization of 
$\alpha(n;k_1,\ldots,k_n)$} Finally, we need the following result from \cite{fischertrapezoids}: let $c,d$ be non-negative integers with $c+d \le n$. We consider the following 
polynomial expansion,
\begin{multline}
\label{expansion}
\alpha(n;k_{1},\ldots,k_{c},c+1,c+2,\ldots,n-d,k_{n-d+1},k_{n-d+2},\ldots,k_{n})
\\
= \sum_{s_{1}=1}^{n} \sum_{s_{2}=1}^{n} \ldots \sum_{s_{c}=1}^{n}
\sum_{i_{1}=1}^{n} \sum_{i_{2}=1}^{n} \ldots \sum_{i_{d}=1}^{n}
B(n;s_{1},s_{2},\ldots,s_{c};i_{1},\ldots,i_{d}) \\
\times (-1)^{s_{1}+s_{2}+\ldots+s_{c}+c}
\binom{k_{1}-c-1}{s_{c}-1} \binom{k_{2}-c-1}{s_{c-1}-1} \cdots
\binom{k_{c}-c-1}{s_{1}-1} \\ \times
\binom{k_{n-d+1}-n+d-2+i_{1}}{i_{1}-1} \binom{k_{n-d+2}-n+d-2+i_{2}}{i_{2}-1}
\cdots \binom{k_{n}-n+d-2+i_{d}}{i_{d}-1},
\end{multline}
where $B(n;s_1,s_2,\ldots,s_c;i_1,\ldots,i_d)$ are certain coefficients. (Note that the 
binomial coefficient $\binom{x}{n} = x (x-1) \cdots (x-n+1) / n!$ is indeed a polynomial in 
$x$ for fixed $n$ and, as $n$ varies, constitutes a basis for the polynomials in $x$.)
Theorem~1 of \cite{fischertrapezoids} states that, for $1 \le s_1 < s_2 < \ldots < s_c \le n$ and 
$1 \le i_1 < i_2 < \ldots < i_d \le n$,  where $B(n;s_1,\ldots,s_c;i_1,\ldots,i_d)=A(n;s_1,\ldots,s_c;i_1,\ldots,i_d)$. 
That is to say that the coefficients enumerate $(n-c-d) \times n$ matrices with entries in $\{1,-1,0\}$ such that the following conditions are fulfilled.
\begin{itemize}
\item The non-zero elements alternate in each  row and column.
\item All rowsums are $1$.
\item The sum of entries in column $j$ is $1$ if and only if $j \notin 
\{i_1,\ldots,i_d,s_1,\ldots,s_c\}$; it is $-1$ if and only if 
$j \in \{i_1,\ldots,i_d\} \cap \{s_1,\ldots,s_c\}$. (In all other 
cases it is $0$.) 
\item The first non-zero entry (if it exists) in columns $i_1,i_2,\ldots,i_d$ is $-1$.
\item The last non-zero entry (if it exists) in columns $s_1, s_2, \ldots, s_c$ is $-1$.
\end{itemize}

\medskip

We also need the following result \cite{fischertrapezoids}, which was in fact used to prove the result that we have just stated (be aware of the fact 
that we interchange the role of rows and columns): 
as usual, let $(s_1,\ldots,s_c)$ and $(i_1,\ldots,i_d)$ be strictly increasing sequences of integers in $\{1,2,\ldots,n\}$, 
$m \ge 1$ and $(k_{c+1}, k_{c+2}, \ldots, k_{n-d})$ a strictly increasing sequence of integers in $\{1,2,\ldots, m\}$.
Then the quantity 
\begin{multline}
\label{lefttopright}
\left( \Delta_{k_1}^{s_c-1} \Delta_{k_2}^{s_{c-1}-1} \cdots \Delta_{k_c}^{s_1-1} 
\delta_{k_{n-d+1}}^{i_1-1} \delta_{k_{n-d+2}}^{i_2-1} \cdots 
\delta_{k_n}^{i_d-1} (-1)^{s_1+\ldots+s_c+c}  \right. \\
\left. \left. \times \alpha(n;k_1,\ldots,k_n)  \right)
\right|_{k_1=k_2=\ldots=k_c=1,k_{n-d+1}=k_{n-d+2}=\ldots=k_n=m}
\end{multline}
(where $\delta_x = \id - E_x^{-1}$) is the number of $n \times m$ matrices with entries in $\{0,1,-1\}$ such
that the non-zero elements in rows and columns alternate and the 
following conditions are fulfilled.
\begin{itemize}
\item The column sum is $1$ precisely for the columns $k_{c+1},k_{c+2},\ldots,
k_{n-d}$; for all other columns the sum is $0$ and the first non-zero entry in these columns (if it exists) is a $-1$.
\item The sum of entries in row $j$ is $1$ if and only if $j \notin 
\{i_1,\ldots,i_d,s_1,\ldots,s_c\}$; it is $-1$ if and only if 
$j \in \{i_1,\ldots,i_d\} \cap \{s_1,\ldots,s_c\}$. (In all other 
cases it is $0$.) 
\item The first non-zero entry (if it exists) in rows $s_1,s_2,\ldots,s_c$ is $-1$.
\item The last non-zero entry (if it exists) in rows $i_1, i_2, \ldots, i_d$ is $-1$.
\end{itemize}

\section{A warm-up: a linear relation between the numbers  $\overline{\overline{A}}_{n,i,j}$ and $\overline{A}|_{n,i,j}$.}
\label{warmup}

To begin with and also to illustrate our method, we first derive another linear relation between two types of 
doubly refined enumerations of alternating sign matrices. (It is in fact also a corollary of Theorem~\ref{main}.) 
Here $\overline{A}|_{n,i,j}$ denotes the number of 
$n \times n$ alternating sign matrices with $1$s in position $(1,i)$ and  $(j,n)$. (Clearly, 
$\overline{A}|_{n,i,j}$ and $|\overline{A}_{n,i,j}$ are trivially related via $\overline{A}|_{n,i,j} = 
|\overline{A}_{n,n+1-i,j}$.)

\medskip

The idea of the proof is in fact very easy: it turns out that both numbers
$\overline{\overline{A}}_{n,i,j}$ and $\overline{A}|_{n,i,j}$ appear in the coefficients of $\alpha(n;1,2,\ldots,n-2, k_{n-1}, k_n)$ as a polynomial in $k_{n-1}$ and $k_{n}$ with respect to a certain polynomial basis. For the first numbers, this 
is a consequence of \eqref{expansion}, see also \eqref{direct}. For the latter numbers, this is a consequence of a combination of Lemma~\ref{step3}, Lemma~\ref{j23} and \eqref{lefttopright}. The linear relation follows by comparing the coefficients.

\medskip

We start by applying Lemma~\ref{step3} ($r=n-1$ and $z=k_{n-1}-n+1$) to see that
\begin{multline*}
\alpha(n;1,2,\ldots,n-2,k_{n-1},k_n) \\
= \sum_{p=0}^{\infty} (-1)^p \binom{k_{n-1}-p}{n-1}
\left. \left( e_p(E_{k_1},\ldots,E_{k_{n-2}},E_{k_n}) \alpha(n;k_1,\ldots,k_n) \right) \right|_{(k_{1},\ldots,k_{n-1})=(1,2,\ldots,n-1)}.
\end{multline*}

\subsection{Eliminating $e_p$ -- another useful lemma.} In order to eliminate $e_p$ from this expression, we need the following result.
In this paper we only need the specials cases $j=2,3$. Note that the statement involves the notation $\sum\limits_{(l_1,\ldots,l_{n-1})}^{(k_1,\ldots,k_n)}$, which we have 
silently introduced in \eqref{recursion}.

\begin{lem} 
\label{j23}
For positive integers $p \ge 0$ and $n \ge j \ge 2$, we have 
\begin{multline*}
\left. \left( e_p(E_{k_1},\ldots,E_{k_{n-j}})  \alpha(n;k_1,\ldots,k_n) \right) \right|_{(k_{1},\ldots,k_{n-j+1})=(1,2,\ldots,n-j+1)} \\
= \sum_{(l_{n-j+1},l_{n-j+2}, \ldots,l_{n-1})}^{(n-j+2,k_{n-j+2},\ldots,k_n)} \sum_{i=1}^{n-j+1} \binom{n-j-i+1}{p} \alpha(n-1;1,2,\ldots,i-1,i+1,\ldots,n-j+1,l_{n-j+1}, \ldots, ,l_{n-1}) \\
+ [p=0] \sum_{(l_{n-j+2},l_{n-j+3},\ldots,l_{n-1})}^{(k_{n-j+2},k_{n-j+3},\ldots,k_n)} \alpha(n-1;1,2,\ldots,n-j+1,l_{n-j+2},\ldots,l_{n-1}).
\end{multline*}
\end{lem}

{\it Proof.} The recursion \eqref{recursion} underlying $\alpha(n;k_1,\ldots,k_n)$ yields the case $p=0$. 
We assume $p>0$ and let
$1 \le i_1 < i_2 < \cdots < i_p \le n-j$.
Then
\begin{multline*}
\left. \left( E_{k_{i_1}} E_{k_{i_2}} \cdots E_{k_{i_p}} \alpha(n;k_1,\ldots,k_n) \right) \right|_{(k_1,\ldots,k_{n-j+1})=(1,2,\ldots,n-j+1)} \\
=
\alpha(n;1,2,\ldots,i_{1}-1,i_1+1,i_1+1,\ldots,i_p-1,i_p+1,i_p+1,\ldots,n-j+1,k_{n-j+2},\ldots,k_n) \\
= \sum_{(l_{n-j+1},\ldots,l_{n-1})}^{(n-j+2,k_{n-j+2},\ldots,k_n)} \sum_{1 \le j_1 \le 2 \le j_2 \le  \ldots \le j_{i_1-1} \le i_1 \atop j_1 < j_2 <  \ldots < j_{i_1-1}} \alpha(n-1;j_1,\ldots,j_{i_1-1},i_1+1,i_1+2,\ldots,n-j+1,l_{n-j+1},\ldots,l_{n-1}) \\
= \alpha(n;1,2,\ldots,i_1-1,i_1+1,i_1+1,\ldots,n-j+1,k_{n-j+2},\ldots,k_n),
\end{multline*}
by \eqref{recursion}. Note that the expression does not depend on $i_2, i_3,\ldots,i_p$. Therefore, the left-hand side in the statement of the lemma is equal to
$$
\sum_{i_1=1}^{n-j} \binom{n-j-i_1}{p-1} \alpha(n;1,2,\ldots,i_1-1,i_1+1,i_1+1,\ldots,n-j+1,k_{n-j+2},\ldots,k_n).
$$
Clearly, we have 
\begin{multline*}
 \alpha(n;1,2,\ldots,i_1-1,i_1+1,i_1+1,\ldots,n-j+1,k_{n-j+2},\ldots,k_n) \\
= \sum_{(l_{n-j+1},\ldots,l_{n-1})}^{(n-j+2,k_{n-j+2},\ldots,k_n)} \sum_{i=1}^{i_1} \alpha(n-1;1,2,\ldots,i-1,i+1,\ldots,n-j+1,l_{n-j+1},\ldots,l_{n-1}), 
\end{multline*}
for $i_1 \le n-j+1$.
The assertion follows from
$$
\sum_{i_1=1}^{n-j} \sum_{i=1}^{i_1} \binom{n-j-i_1}{p-1} = \sum_{i=1}^{n-j} \sum_{i_1=i}^{n-j} \binom{n-j-i_1}{p-1} =
\sum_{i=1}^{n-j} \sum_{i_1=0}^{n-j-i} \binom{i_1}{p-1} = 
 \binom{n-j-i+1}{p}. \qed
$$

\medskip

\subsection{Derivation of the linear relation}
For the special case $j=2$, this gives
\begin{multline*}
\left. \left( e_p(E_{k_1},\ldots,E_{k_{n-2}}) \alpha(n;k_1,\ldots,k_n) \right) \right|_{(k_{1},\ldots,k_{n-1})=(1,2,\ldots,n-1)} \\
= \sum_{l_{n-1}=n}^{k_n} \sum_{i=1}^{n-1} \binom{n-1-i}{p} \alpha(n-1;1,2,\ldots,i-1,i+1,\ldots,n-1,l_{n-1})
+ [p=0] A_{n-1}.
\end{multline*}

\medskip


Consequently,
\begin{multline*}
\alpha(n;1,2,\ldots,n-2,k_{n-1},k_n) \\
= \sum_{p=0}^{\infty} (-1)^p \binom{k_{n-1}-p}{n-1}
\left( \left.  \left( e_p(E_{k_1},\ldots,E_{k_{n-2}}) \alpha(n;k_1,\ldots,k_n) \right) \right|_{(k_{1},\ldots,k_{n-1})=(1,2,\ldots,n-1)} \right. \\
\left. \left. + \left( e_{p-1}(E_{k_1},\ldots,E_{k_{n-2}}) \alpha(n;k_1,\ldots,k_{n-1},k_n+1) \right) \right|_{(k_{1},\ldots,k_{n-1})=(1,2,\ldots,n-1)} \right) \\
= \sum_{p=0}^{\infty} (-1)^p \binom{k_{n-1}-p}{n-1} \sum_{l_{n-1}=n}^{k_n} \sum_{i=1}^{n-1}
\binom{n-1-i}{p} \alpha(n-1;1,\ldots,i-1,i+1,\ldots,n-1,l_{n-1}) \\ + A_{n-1} \binom{k_{n-1}}{n-1}  \\
+ \sum_{p=1}^{\infty} (-1)^p \binom{k_{n-1}-p}{n-1} \sum_{l_{n-1}=n}^{k_n+1} \sum_{i=1}^{n-1} \binom{n-1-i}{p-1}
\alpha(n-1;1,\ldots,i-1,i+1,\ldots,n-1,l_{n-1}) \\ -  A_{n-1} \binom{k_{n-1}-1}{n-1}.
\end{multline*}

We use the following identity, which can be deduced from the Chu-Vandermonde summation and holds for integers $n,i,j$ with $i \le n-1$ and the indeterminate $k_{n-1}$,
$$
\sum_{p=j}^{\infty} (-1)^p \binom{k_{n-1}-p}{n-1} \binom{n-1-i}{p-j} =
(-1)^j \binom{k_{n-1} -n +1 - j +i}{i}
$$
to see that this is equal to
\begin{multline*}
\sum_{l_{n-1}=n}^{k_n} \sum_{i=1}^{n-1} \binom{k_{n-1}-n+1+i}{i} \alpha(n-1;1,\ldots,i-1,i+1,\ldots,n-1,l_{n-1}) \\
- \sum_{l_{n-1}=n}^{k_n+1} \sum_{i=1}^{n-1} \binom{k_{n-1}-n+i}{i} \alpha(n-1;1,\ldots,i-1,i+1,\ldots,n-1,l_{n-1}) \\
+ \binom{k_{n-1}-1}{n-2} A_{n-1}.
\end{multline*}
Clearly, $\alpha(n-1;1,\ldots,i-1,i+1,\ldots,n-1,l_{n-1})$ is a polynomial in $l_{n-1}$ of degree no greater than
$n-2$ and therefore possesses an expansion in terms of the polynomial basis $\left( \binom{l_{n-1}-n-1+j}{j-1} \right)_{1 \le j \le n-1}$.
$$
\alpha(n-1;1,\ldots,i-1,i+1,\ldots,n-1,l_{n-1}) = \sum_{j=1}^{n} B_{n,i,j} \binom{l_{n-1}-n-2+j}{j-2}
$$
However, the  coefficient $B_{n,i,j}$  is the number of $n \times n$ alternating sign matrices where the unique $1$ in the first row is in
column $i$ and the unique $1$ in the last column is in row $j$, i.e. $B_{n,i,j}= \overline{A}|_{n,i,j}$. This follows from
$$
\overline{A}|_{n,i,j} = \left. \delta_{l_{n-1}}^{j-2} \alpha(n-1;1,\ldots,i-1,i+1,\ldots,n-1,l_{n-1}) \right|_{l_{n-1}=n-1},
$$
(where $\delta_x = \id - E^{-1}_x$) which is a special case of \eqref{lefttopright}.

\medskip

Therefore,
\begin{multline*}
\alpha(n;1,2,\ldots,n-2,k_{n-1},k_n) \\
= \sum_{l_{n-1}=n}^{k_n} \sum_{i=1}^{n-1} \binom{k_{n-1}-n+1+i}{i} \sum_{j=1}^{n} \overline{A}|_{n,i,j} \binom{l_{n-1}-n-2+j}{j-2} \\
- \sum_{l_{n-1}=n}^{k_n+1} \sum_{i=1}^{n-1} \binom{k_{n-1}-n+i}{i} \sum_{j=1}^{n} \overline{A}|_{n,i,j}
\binom{l_{n-1}-n-2+j}{j-2} + A_{n-1} \binom{k_{n-1}-1}{n-2}  \\
=  \sum_{i=2}^{n} \sum_{j=1}^{n} \overline{A}|_{n,i-1,j} \binom{k_{n-1}-n+i}{i-1} \binom{k_n -n + j-1}{j-1} \\
- \sum_{i=2}^{n} \sum_{j=1}^{n} \overline{A}|_{n,i-1,j} \binom{k_{n-1} - n + i-1}{i-1} \binom{k_n - n + j}{j-1} \\
 + A_{n-1} \binom{k_{n-1}-1}{n-2},
\end{multline*}
where we have used the fact that 
$$
\sum_{x=a}^{b} \binom{x}{n} = \sum_{x=a}^{b} \left( \binom{x+1}{n+1} - \binom{x}{n+1} \right) = \binom{b+1}{n+1} - \binom{a}{n+1}.
$$
We use
$$
\binom{k_{n-1}-n+i-1}{i-1} = \binom{k_{n-1}-n+i}{i-1} - \binom{k_{n-1}-n+i-1}{i-2}
$$
and
$$
\binom{k_n-n+j-1}{j-1} = \binom{k_n-n+j}{j-1} - \binom{k_n-n+j-1}{j-2}, 
$$
to see that this implies
\begin{multline*}
\alpha(n;1,2,\ldots,n-2,k_{n-1},k_n)
=
- \sum_{i=2}^{n} \sum_{j=0}^{n-1} \overline{A}|_{n,i-1,j+1} \binom{k_{n-1}-n+i}{i-1} \binom{k_n -n + j}{j-1} \\
+ \sum_{i=1}^{n-1} \sum_{j=1}^{n} \overline{A}|_{n,i,j} \binom{k_{n-1} - n + i}{i-1} \binom{k_n - n + j}{j-1}
+ A_{n-1} \binom{k_{n-1}-1}{n-2}
\end{multline*}
On the other hand, by \eqref{expansion}, 
\begin{equation}
\label{direct}
\alpha(n;1,2,\ldots,n-2,k_{n-1},k_n) = \sum_{i=1}^{n} \sum_{j=1}^{n} \overline{\overline{A}}_{n,i,j} \binom{k_{n-1}-n+i}{i-1}
\binom{k_n-n+j}{j-1}.
\end{equation}
If we compare the coefficients, this yields
\begin{equation*}
\overline{\overline{A}}_{n,i,j} = -\overline{A}|_{n,i-1,j+1} + \overline{A}|_{n,i,j} + ([i=n-1][j=1] - [i=n] [j=1]) A_{n-1}
\end{equation*}
for $i, j \in \{1,2,\ldots,n\}$ where we set $\overline{A}|_{n,i,j}=0$ if $i=0$ or $j=n+1$. Observe that this identity follows from Theorem~\ref{main} if we set $s=n$. Furthermore: if we set $\overline{\overline{A}}_{n,i,j}=0$ and $\overline{\overline{A}}_{n,i,j}=0$ whenever $i \notin \{1,2,\ldots,n\}$ or $j \notin \{1,2,\ldots,n\}$ then the following modified identity holds true for all integers $i$ and $j$.
\begin{equation}
\label{doublyilse}
\overline{\overline{A}}_{n,i,j} = (\id - E_j E^{-1}_i) \overline{A}|_{n,i,j}  + ([i=n-1][j=1] - [i=n] [j=1] + [i=n+1] [j=0]) A_{n-1}
\end{equation}

\medskip

\subsection{Reversing the linear relations} Of course, this identity can also be used to derive a formula for $\overline{A}|_{n,i,j}$ in terms of the numbers
$\overline{\overline{A}}_{n,i,j}$. We have 
\begin{multline*}
\overline{A}|_{n,i,j} = (\id - E_j E^{-1}_i)^{-1} \overline{\overline{A}}_{n,i,j}
\\ + (\id - E_j E^{-1}_i)^{-1} ( - [i=n-1] [j=1] + [i=n] [j=1] - [i=n+1] [j=0]) 
A_{n-1} \\
= \sum_{k=0}^{\infty}  \overline{\overline{A}}_{n,i-k,j+k} +
(- [i+j=n] [j \le 1] + [i=n] [j=1]) A_{n-1}, 
\end{multline*}
where the sum is in fact finite as  $\overline{\overline{A}}_{n,i-k,j+k}=0$ if 
$k \ge \max(i+1,n-j+1)$, which is also the reason why we may use the 
expansion $(\id - E_j E^{-1}_i)^{-1} = \sum\limits_{k=0}^{\infty} E^k_j E^{-k}_i$.
To more precise, the following proposition lurks behind.

\begin{prop}
\label{elegant1}
 Let $n$ be an integer.
Suppose that $(B_{i,j})_{i, j \in \mathbb{Z}}$ is an infinite matrix with 
$B_{i,j} = 0$ if $i \le 0$ or $j \ge n+1$. Then there is a unique infinite
matrix $(A_{i,j})_{i, j \in \mathbb{Z}}$ with 
$$
B_{i,j} = (\id - E_j E^{-1}_i) A_{i,j} 
$$
and $A_{0,j} = A_{i,n+1} = 0$ for all $i, j \in \mathbb{Z}$. It is given by
$$
A_{i,j} = \sum_{k=0}^{\infty} B_{i-k,j+k}.
$$
\end{prop}

\medskip

Another possibility to derive a formula for 
$\overline{A}|_{n,i,j}$ is the following. We apply $E^{-1}_j E_i$ to both sides of \eqref{doublyilse} and rearrange the identity to arrive at
\begin{multline*}
(\id - E_j^{-1} E_i) \overline{A}|_{n,i,j} = - E_j^{-1} E_i 
\overline{\overline{A}}_{n,i,j} \\
+ ( [i+1=n-1] [j-1=1]  - [i+1=n] [j-1 = 1] + [i+1=n+1] [j-1=0]) A_{n-1}.
\end{multline*}
We apply $(\id - E_j^{-1} E_i)^{-1}$ to see that 
\begin{multline*}
\overline{A}|_{n,i,j} = - \sum_{k=0}^{\infty} 
\overline{\overline{A}}_{n,i+1+k,j-1-k} \\ 
+ \sum_{k=0}^{\infty} 
([i+k+1=n-1] [j-k-1=1]  - [i+k+1=n] [j-k-1=1] + [i+k+1=n+1] [j-k-1=0]) A_{n-1}
\\
=  - \sum_{k=0}^{\infty} 
\overline{\overline{A}}_{n,i+1+k,j-1-k} + ([i+j=n] [j \ge 2] + [i=n] [j=1]) A_{n-1}
\end{multline*}
Interestingly, these two formulas for $\overline{A}|_{n,i,j}$ imply the following identity for the numbers $\overline{\overline{A}}_{n,i,j}$.
$$
\sum_{i,j \in \mathbb{Z} \atop i+j=p} \overline{\overline{A}}_{n,i,j} = 
[p=n] A_{n-1}.
$$

\medskip

Finally, we use \eqref{stroganov} and \eqref{doublyilse} to derive a direct relation between $\overline{\overline{A}}_{n,i,j}$ and $\overline{\underline{A}}_{n,i,j}$ of the same type as \eqref{kark}. In \eqref{doublyilse}, we replace 
$i$ by $n+1-i$ and use the fact that $\overline{A}|_{n,n+1-i,j} = 
|\overline{A}_{n,i,j}$ to obtain 
\begin{equation}
\label{first}
\overline{\overline{A}}_{n,n+1-i,j} = (\id - E_j E_i) 
|\overline{A}_{n,i,j} + ([i=0] [j=0] - [i=1] [j=1] + [i=2] [j=1]) A_{n-1}
\end{equation}
for all integers $i$ and $j$. 
We may extend
\eqref{stroganov} to all integers $i, j$ by modifying it as follows
\begin{equation}
\label{second}
\overline{\underline{A}}_{n,i,j} = 
(E_i + E_j - E_i E_j) |\overline{A}_{n,i,j} + ([i=1] - [i=0]) ( [j=1] - [j=0]) A_{n-1}.
\end{equation}
We apply the operator $(E_i + E_j - E_i E_j)$ to \eqref{first} and use \eqref{second} to express $(E_i + E_j - E_i E_j) |\overline{A}_{n,i,j}$ in terms of 
$\overline{\underline{A}}_{n,i,j}$.
This yields
$$
(E_i + E_j - E_i E_j) \overline{\overline{A}}_{n,n+1-i,j} = 
(\id - E_j E_i) \overline{\underline{A}}_{n,i,j} + ([i=2] [j=0] + [i=1] [j=0]
+ [i=0] [j=0]) A_{n-1}.
$$
(Note that in contrary to \eqref{kark}, we do not necessarily have  $\overline{\overline{A}}_{n,i,j} \not= 0$ if $i > j$.)

\medskip 

Finally, \eqref{second} can also be used to express $|\overline{A}_{n,i,j}$ in 
terms of $\overline{\underline{A}}_{n,i,j}$. We apply $E^{-1}_i E^{-1}_j$ to both sides of the equation and multiply the result by $-1$. In order to apply 
$(\id - (E_i^{-1} + E_j^{-1})^{-1}$, we need a proposition that is analogous to Proposition~\ref{elegant1}:
\begin{prop}
\label{elegant}
Suppose that $(B_{i,j})_{i, j \in \mathbb{Z}}$ is an infinite matrix with 
$B_{i,j} = 0$ if $i \le 0$ or $j \le 0$. Then there is a unique infinite
matrix $(A_{i,j})_{i, j \in \mathbb{Z}}$ with 
$$
B_{i,j} = (\id - (E_i^{-1} + E^{-1}_j)) A_{i,j} 
$$
and $A_{0,j} = A_{i,0} = 0$ for all $i, j \in \mathbb{Z}$. It is given by
$$
A_{i,j} =\sum_{k=0}^{\infty} \sum_{l=0}^{k} \binom{k}{l}
 B_{i-l,j-k+l}.
$$
\end{prop}
This  leads to 
$$
|\overline{A}_{n,i,j} = - \sum_{k=0}^{\infty} \sum_{l=0}^{k} \binom{k}{l} 
\overline{\underline{A}}_{n,i-l-1,j-k+l-1} + A_{n-1} \binom{i+j-4}{i-2} [i \ge 2] [j \ge 2] + 
A_{n-1} [i=1] [j=1].
$$

\section{Proof of Theorem~\ref{main}}
\label{triply}

The proof of Theorem~\ref{main} is very similar to the proof of \eqref{doublyilse}. However, the computations are more involved, although most steps are really elementary. 
In this case
\eqref{cycle} and \eqref{expansion} imply easily that the numbers 
$\overline{\overline{\underline{A_{n,s,i,t}}}}$ appear as the coefficient 
of $\alpha(n;1,2,\ldots,n-3,k_{n-2},k_{n-1},k_n)$ with respect to a certain 
polynomial basis. The fact that also $|\overline{A}|_{n,s,i,t}$ appears in the coefficient of $\alpha(n;1,2,\ldots,n-3,k_{n-2},k_{n-1},k_n)$ is again a consequence of Lemma~\ref{step3}, Lemma~\ref{j23} (now we need the case $j=3$) and \eqref{lefttopright}. 

\medskip

\subsection{Applying Lemma~\ref{step3} -- introducing $e_p$}
Lemma~\ref{step3} ($r=n-2$ and $z=k_{n-2}-n+2$) implies
\begin{multline*}
\alpha(n;1,2,\ldots,n-3,k_{n-2},k_{n-1},k_n) \\
= \sum_{p=0}^{\infty} (-1)^p \binom{k_{n-2}-p+1}{n-1}
\left. \left( e_p(E_{k_1},\ldots,E_{k_{n-3}},E_{k_{n-1}},E_{k_n}) \alpha(n;k_1,\ldots,k_n) \right) \right|_{(k_{1},\ldots,k_{n-2})=(1,2,\ldots,n-2)}.
\end{multline*}
This is  equal to
\begin{multline*}
\alpha(n;1,2,\ldots,n-3,k_{n-2},k_{n-1},k_n) \\
= \sum_{p=0}^{\infty} (-1)^p \binom{k_{n-2}-p+1}{n-1}
\left. \left( e_p(E_{k_1},\ldots,E_{k_{n-3}}) \alpha(n;k_1,\ldots,k_n) \right) \right|_{(k_{1},\ldots,k_{n-2})=(1,2,\ldots,n-2)} \\
+ \sum_{p=1}^{\infty} (-1)^p \binom{k_{n-2}-p+1}{n-1}
\left. \left( e_{p-1}(E_{k_1},\ldots,E_{k_{n-3}}) \alpha(n;k_1,\ldots,k_{n-1}+1,k_n) \right) \right|_{(k_{1},\ldots,k_{n-2})=(1,2,\ldots,n-2)} \\
+ \sum_{p=1}^{\infty} (-1)^p \binom{k_{n-2}-p+1}{n-1}
\left. \left( e_{p-1}(E_{k_1},\ldots,E_{k_{n-3}}) \alpha(n;k_1,\ldots,k_{n-1},k_n+1) \right) \right|_{(k_{1},\ldots,k_{n-2})=(1,2,\ldots,n-2)} \\
+ \sum_{p=2}^{\infty} (-1)^p \binom{k_{n-2}-p+1}{n-1}
\left. \left( e_{p-2}(E_{k_1},\ldots,E_{k_{n-3}}) \alpha(n;k_1,\ldots,k_{n-1}+1,k_n+1) \right) \right|_{(k_{1},\ldots,k_{n-2})=(1,2,\ldots,n-2)}.
\end{multline*}

\medskip

\subsection{Applying Lemma~\ref{j23} -- eliminating $e_p$} 
Again, we use Lemma~\ref{j23} to eliminate $e_p$ on 
the right-hand side of this identity: if we set $j=3$ in the lemma then we 
obtain
\begin{multline*}
\left. \left( e_p(E_{k_1},\ldots,E_{k_{n-3}}) \alpha(n;k_1,\ldots,k_n) \right) \right|_{(k_{1},\ldots,k_{n-2})=(1,2,\ldots,n-2)} \\
= \sum_{(l_{n-2},l_{n-1})}^{(n-1,k_{n-1},k_n)} \sum_{i=1}^{n-2} \binom{n-2-i}{p} \alpha(n-1;1,2,\ldots,i-1,i+1,\ldots,n-2,l_{n-2},l_{n-1}) \\
+ [p=0] \sum_{l_{n-1}=k_{n-1}}^{k_n} \alpha(n-1;1,2,\ldots,n-2,l_{n-1}).
\end{multline*}
Therefore and after some further simplifications, we obtain
\begin{multline*}
\alpha(n;1,2,\ldots,n-3,k_{n-2},k_{n-1},k_n) \\
= \sum_{p=0}^{\infty} (-1)^p \binom{k_{n-2}-p+1}{n-1}
  \sum_{(l_{n-2},l_{n-1})}^{(n-1,k_{n-1},k_n)} \sum_{i=1}^{n-2} \binom{n-2-i}{p} \alpha(n-1;1,\ldots,i-1,i+1,\ldots,n-2,l_{n-2},l_{n-1})  \\
+ \sum_{p=1}^{\infty} (-1)^p \binom{k_{n-2}-p+1}{n-1}
\sum_{(l_{n-2},l_{n-1})}^{(n-1,k_{n-1}+1,k_n)} \sum_{i=1}^{n-2} \binom{n-2-i}{p-1} \alpha(n-1;1,\ldots,i-1,i+1,\ldots,n-2,l_{n-2},l_{n-1})  \\
+ \sum_{p=1}^{\infty} (-1)^p \binom{k_{n-2}-p+1}{n-1}
\sum_{(l_{n-2},l_{n-1})}^{(n-1,k_{n-1},k_n+1)} \sum_{i=1}^{n-2} \binom{n-2-i}{p-1} \alpha(n-1;1,\ldots,i-1,i+1,\ldots,n-2,l_{n-2},l_{n-1})  \\
+ \sum_{p=2}^{\infty} (-1)^p \binom{k_{n-2}-p+1}{n-1}
\sum_{(l_{n-2},l_{n-1})}^{(n-1,k_{n-1}+1,k_n+1)} \sum_{i=1}^{n-2} \binom{n-2-i}{p-2} \alpha(n-1;1,\ldots,i-1,i+1,\ldots,n-2,l_{n-2},l_{n-1})  \\
+ \binom{k_{n-2}-1}{n-3} \sum_{l_{n-1}=k_{n-1}}^{k_n} \alpha(n-1;1,2,\ldots,n-2,l_{n-1}) \\
+ \binom{k_{n-2}-1}{n-2} (\alpha(n-1;1,2,\ldots,n-2,k_{n-1}) - \alpha(n-1;1,2,\ldots,n-2,k_{n}+1))
\end{multline*}
Now  we need the following identity, which follows from the Chu--Vandermonde summation
and holds all integers $n, i, j$ with $i \le n-2$ and the indeterminate $k_{n-2}$,
$$
\sum_{p=j}^{\infty} (-1)^{p} \binom{k_{n-2}-p+1}{n-1} \binom{n-2-i}{p-j} =
(-1)^{j} \binom{k_{n-2}-n+i-j+3}{i+1}
$$
to see that
\begin{multline}
\label{shift}
\alpha(n;1,2,\ldots,n-3,k_{n-2},k_{n-1},k_n) \\
=  \sum_{(l_{n-2},l_{n-1})}^{(n-1,k_{n-1},k_n)} \sum_{i=2}^{n-1} \binom{k_{n-2}-n+i+2}{i} \alpha(n-1;1,2,\ldots,i-2,i,\ldots,n-2,l_{n-2},l_{n-1})  \\
- \sum_{(l_{n-2},l_{n-1})}^{(n-1,k_{n-1}+1,k_n)} \sum_{i=2}^{n-1} \binom{k_{n-2}-n+i+1}{i} \alpha(n-1;1,2,\ldots,i-2,i,\ldots,n-2,l_{n-2},l_{n-1})  \\
- \sum_{(l_{n-2},l_{n-1})}^{(n-1,k_{n-1},k_n+1)} \sum_{i=2}^{n-1} \binom{k_{n-2}-n+i+1}{i} \alpha(n-1;1,2,\ldots,i-2,i,\ldots,n-2,l_{n-2},l_{n-1})  \\
+ \sum_{(l_{n-2},l_{n-1})}^{(n-1,k_{n-1}+1,k_n+1)} \sum_{i=2}^{n-1} \binom{k_{n-2}-n+i}{i} \alpha(n-1;1,2,\ldots,i-2,i,\ldots,n-2,l_{n-2},l_{n-1})  \\
+ \binom{k_{n-2}-1}{n-3} \sum_{l_{n-1}=k_{n-1}}^{k_n} \alpha(n-1;1,2,\ldots,n-2,l_{n-1}) \\
+ \binom{k_{n-2}-1}{n-2} (\alpha(n-1;1,2,\ldots,n-2,k_{n-1}) - \alpha(n-1;1,2,\ldots,n-2,k_{n}+1)).
\end{multline}
(For what follows a shift of $i$ is useful.)

\medskip

\subsection{The extra term $R$} Next we deal with the extra term in \eqref{shift}, i.e. with
\begin{multline}
\label{rest}
\binom{k_{n-2}-1}{n-3} \sum_{l_{n-1}=k_{n-1}}^{k_n} \alpha(n-1;1,2,\ldots,n-2,l_{n-1}) \\
+ \binom{k_{n-2}-1}{n-2} (\alpha(n-1;1,2,\ldots,n-2,k_{n-1}) -
\alpha(n-1;1,2,\ldots,n-2,k_{n}+1))=:R,
\end{multline}
which we denote by $R$ in the following. By \eqref{expansion}, we have 
$$
\alpha(n-1;1,2,\ldots,n-2,x) = \sum_{l=1}^{n-1} A_{n-1,l} \binom{x-n+l}{l-1}.
$$
Thus, \eqref{rest} is equal to
\begin{multline*}
\binom{k_{n-2}-1}{n-3} \sum_{l=1}^{n-1} A_{n-1,l}
\sum_{l_{n-1}=k_{n-1}}^{k_n}  \binom{l_{n-1}-n+l}{l-1} \\
+ \binom{k_{n-2}-1}{n-2} \sum_{l=1}^{n-1} A_{n-1,l} \left( \binom{k_{n-1}-n+l}{l-1}
- \binom{k_n-n+l+1}{l-1} \right) = \\
\binom{k_{n-2}-1}{n-3} \sum_{l=1}^{n-1} A_{n-1,l} \left( \binom{k_n-n+l+1}{l} -
\binom{k_{n-1} - n +l}{l} \right) \\
+ \binom{k_{n-2}-1}{n-2} \sum_{l=1}^{n-1} A_{n-1,l} \left( \binom{k_{n-1}-n+l}{l-1}
- \binom{k_n-n+l+1}{l-1} \right)
\end{multline*}
Using this, and standard identities for binomial coefficients, it is not hard to see that 
\begin{multline*}
R=  \binom{k_{n-2}-1}{n-3} \sum_{t=1}^{n} (A_{n-1,t} - A_{n-1,t-1}) \binom{k_{n-1}-n+t+1}{t-1} \\
+ \binom{k_{n-2}}{n-2} \sum_{t=1}^{n} (A_{n-1,t} - A_{n-1,t+1}) \binom{k_{n-1}-n+t+1}{t-1} \\
+ \binom{k_{n-2}-1}{n-3} \sum_{j=1}^{n} A_{n-1,j-1} \binom{k_{n} - n + j+1}{j-1}
- \binom{k_{n-2}}{n-2} \sum_{j=1}^{n} A_{n-1,j} \binom{k_n-n+j+1}{j-1},
\end{multline*}
where $A_{n-1,0} = A_{n-1,n} = 0$.
In order to modify $R$ a little further, we use 
$$
\binom{k_n-n+j+1}{j-1} = \sum_{s=1}^{j} \binom{-2n+2+j}{j-s} \binom{k_n+n-1}{s-1}
$$
(which follows from the Chu-Vandermonde identity), to see that 
\begin{multline*}
\sum_{j=1}^{n} A_{n-1,j-1} \binom{k_n-n+j+1}{j-1} =
\sum_{s=1}^{n} \sum_{j=s}^{n} A_{n-1,j-1} \binom{-2n+2+j}{j-s} \binom{k_n+n-1}{s-1} \\
= \sum_{s=1}^{n} \sum_{t=s-1}^{n-1} A_{n-1,t} \binom{2n-3-s}{t+1-s} (-1)^{s+t+1} 
\binom{k_n+n-1}{s-1} \\
= \sum_{s=1}^{n} \sum_{t=s-1}^{n-1} A_{n-1,t} \left( \binom{2n-2-s}{t+1-s} - \binom{2n-3-s}{t-s} \right) (-1)^{s+t+1} \binom{k_n+n-1}{s-1}.
\end{multline*}
Now we need the following identity, which is a special case of 
\eqref{topbottom1},
\begin{equation}
\label{single}
A_{n-1,s} = \sum_{t=s}^{n-1} (-1)^{n+t+1} A_{n-1,t}
\binom{2n-3-s}{t-s}
\end{equation}
to see that this is furthermore equal to
$$
\sum_{s=1}^{n} (-1)^{n+s} (A_{n-1,s-1} - A_{n-1,s} ) \binom{k_n+n-1}{s-1}.
$$
Similarly, 
\begin{multline*}
\sum_{j=1}^{n} A_{n-1,j} \binom{k_n-n+j+1}{j-1} 
= \sum_{s=1}^{n} \sum_{j=s}^{n} A_{n-1,j} \binom{-2n+2+j}{j-s} \binom{k_n+n-1}{s-1} \\
= \sum_{s=1}^{n} \sum_{j=s}^{n-1} A_{n-1,j} \binom{2n-3-s}{j-s} (-1)^{j+s} 
\binom{k_n+n-1}{s-1} = \sum_{s=1}^{n} A_{n-1,s} (-1)^{n+1+s} \binom{k_n+n-1}{s-1}.
\end{multline*}
Therefore 
\begin{multline*}
R = \binom{k_{n-2}-1}{n-3} \sum_{t=1}^{n} (A_{n-1,t} - A_{n-1,t-1}) \binom{k_{n-1}-n+t+1}{t-1} \\
+ \binom{k_{n-2}}{n-2} \sum_{t=1}^{n} (A_{n-1,t} - A_{n-1,t+1}) \binom{k_{n-1}-n+t+1}{t-1} \\
+ \binom{k_{n-2}-1}{n-3} \sum_{s=1}^ {n} (-1)^{n+s} (A_{n-1,s-1} - A_{n-1,s}) \binom{k_n+n-1}{s-1} \\
+ \binom{k_{n-2}}{n-2} \sum_{s=1}^{n} (-1)^{n+s} A_{n-1,s} \binom{k_n+n-1}{s-1}.
\end{multline*}

\medskip

\subsection{More facts about $\alpha(n;k_1,\ldots,k_n)$}
Using \eqref{operator} as well as the definition of $\alpha(n;k_1,\ldots,k_n)$, we see that
\begin{multline*}
\alpha(n;k_1,\ldots,k_n) = \alpha(n;k_1-n+1,k_2-n+1,\ldots,k_n-n+1) = \\
\prod_{1 \le p < q \le n} E^{-1}_{k_p} E^{-1}_{k_q} \prod_{1 \le p < q \le n} (\id - E_{k_p} + E_{k_p} E_{k_q}) \prod_{1 \le i < j \le n}
\frac{k_j - k_i}{j-i} = \\
\prod_{1 \le p < q \le n} (\id - E^{-1}_{k_q} + E^{-1}_{k_p} E^{-1}_{k_q}) \prod_{1 \le i < j \le n}
\frac{k_j - k_i}{j-i} =
\prod_{1 \le p < q \le n} (\id - \delta_{k_p} + \delta_{k_p} \delta_{k_q}) \prod_{1 \le i < j \le n}
\frac{k_j - k_i}{j-i}.
\end{multline*}
As
$$
\prod_{1 \le i < j \le n}
\frac{k_j - k_i}{j-i} = \det_{1 \le i, j \le n}  \binom{k_i -n + j}{j-1}
$$
we also have
$$
\alpha(n;k_1,\ldots,k_n) = \prod_{1 \le p < q \le n} (\id - \delta_{k_p} + \delta_{k_p} \delta_{k_q})
\det_{1 \le i, j \le n}  \binom{k_i -n + j}{j-1}.
$$
For $m \ge 0$, we set
\begin{equation}
\label{m}
\alpha_m(n;k_1,\ldots,k_n) = \prod_{1 \le p < q \le n} (\id - \delta_{k_p} + \delta_{k_p} \delta_{k_q})
\det_{1 \le i, j \le n}  \binom{k_i -n + j+m}{j+m-1}.
\end{equation}
Note that $\alpha_0(n;k_1,\ldots,k_n) = \alpha(n;k_1,\ldots,k_n)$ and
$$
\delta_{k_1} \delta_{k_2} \cdots \delta_{k_n} \alpha_m(n;k_1,\ldots,k_n) =
\alpha_{m-1}(n;k_1,\ldots,k_n).
$$
Therefore and as \eqref{m} can easily be used to obtain the expansion of $\alpha_m(n;k_1,\ldots,k_n)$ in
terms of the polynomial basis
$$\binom{k_1-n+j_1}{j_1-1} \binom{k_2-n+j_2}{j_2-1} \cdots \binom{k_n-n+j_n}{j_n-1} \qquad (j_1,j_2,\ldots, j_n \ge 0),$$
(to this end note also that $\delta_{k_i} \binom{k_i-n+j}{j-1} = \binom{k_i-n+j-1}{j-2}$) it is obvious that $\alpha_{m+1}(n;k_1,\ldots,k_n)$ can be obtained from this expansion
of $\alpha_{m}(n;k_1,\ldots,k_n)$ by replacing
$$\binom{k_1-n+j_1}{j_1-1} \binom{k_2-n+j_2}{j_2-1} \cdots \binom{k_n-n+j_n}{j_n-1}$$
with
$$\binom{k_1-n+j_1+1}{j_1} \binom{k_2-n+j_2+1}{j_2} \cdots \binom{k_n-n+j_n+1}{j_n}.$$
More general:
$\delta_{k_{i_1}} \delta_{k_{i_2}} \cdots \delta_{k_{i_d}} \alpha_{m+1}(n;k_1,\ldots,k_n)$ can be obtained from $\alpha_m(n;k_1,\ldots,k_n)$ by replacing
$$\binom{k_1-n+j_1}{j_1-1} \binom{k_2-n+j_2}{j_2-1} \cdots \binom{k_n-n+j_n}{j_n-1}$$
with
$$\prod_{l=1}^{n} \binom{k_l-n+j_l+1}{j_l}
\prod_{l=1}^{d} \frac{\binom{k_{i_l}-n+j_{i_l}}{j_{i_l}-1}}
{\binom{k_{i_l}-n+j_{i_l}+1}{j_{i_l}}}.$$

\medskip

This all is useful, because \cite[Lemma~1]{fischergtriangle}
implies that
\begin{multline}
\label{useful}
\sum_{(l_i,l_{i+1},\ldots,l_j)}^{(k_i,k_{i+1},\ldots,k_{j+1})}
\alpha_m(n;l_1,\ldots,l_n) =
\sum_{(l_i,l_{i+1},\ldots,l_j)}^{(k_i,k_{i+1},\ldots,k_{j+1})}
\delta_{l_1} \delta_{l_2} \ldots \delta_{l_n}
\alpha_{m+1}(n;l_1,\ldots,l_n) \\
= \delta_{l_1} \ldots \delta_{l_{i-1}} \delta_{l_{j+1}} \ldots \delta_{l_n}
\sum_{r=i}^{j+1} (-1)^{r+i} \alpha_{m+1}(n;l_1,\ldots,l_{i-1},k_i-1,\ldots,k_{r-1}-1,k_{r+1},\ldots,k_{j+1},l_{j+1},\ldots,l_n),
\end{multline}
as, by the definition of
$\alpha_m(n;k_1,\ldots,k_n)$, we know that  $$V_{k_i,k_{i+1}} \alpha_m(n;k_1,\ldots,k_n)$$ is antisymmetric in $k_i$ and $k_{i+1}$, 
where $V_{x,y} = E_x + \Delta_x \Delta_y$.

\medskip

\subsection{The numbers $|\overline{A}|_{n,s,i,t}$ come into the play} 
Equation~\eqref{cycle} implies
\begin{multline*}
\alpha(n-1;1,2,\ldots,i-2,i,\ldots,n-2,l_{n-2},l_{n-1})  \\
= (-1)^{n} \alpha(n-1;l_{n-1}+n-1,1,2,\ldots,i-2,i,\ldots,n-2,l_{n-2}) \\
= (-1)^{n} \alpha(n-1;l_{n-1}+n,2,3,\ldots,i-1,i+1,\ldots,n-1,l_{n-2}+1).
\end{multline*}
Furthermore, by \eqref{lefttopright}, 
$$
|\overline{A}|_{n,s,i,t} = \left. (-1)^s \Delta^{s-2}_{k_1} \delta^{t-2}_{k_{n-1}} \alpha(n-1;k_1,2,\ldots,i-1,i+1,\ldots,n-1,k_{n-1}) \right|_{k_1=2, k_{n-1}=n-1}
$$
and therefore these numbers appear as the coefficient in the following expansion
$$
\alpha(n-1;k_1,2,3,\ldots,i-1,i+1,\ldots,n-1,k_{n-1}) = \sum_{s=1}^{n} \sum_{t=1}^{n}
|\overline{A}|_{n,s,i,t} (-1)^{s} \binom{k_1-2}{s - 2} \binom{k_{n-1}-n-2+t}{t-2}.
$$
We combine this with the previous observation to see that
\begin{multline*}
\alpha(n-1;1,2,\ldots,i-2,i,\ldots,n-2,l_{n-2},l_{n-1}) \\
=  \sum_{s=1}^{n} \sum_{t=1}^{n} |\overline{A}|_{n,s,i,t} (-1)^{n+s} \binom{l_{n-1}+n-2}{s-2} \binom{l_{n-2}-n-1+t}{t-2}.
\end{multline*}
We need to expand $\binom{l_{n-1}+n-2}{s-2}$ in terms of the basis
$(\binom{l_{n-1}-n+r}{r-1})_{r \ge 0}$. This is accomplished by using the Chu--Vandermonde summation.
$$
\binom{l_{n-1}+n-2}{s-2} = \sum_{r=1}^{s-1} (-1)^{r+s+1} \binom{-2n+s+1}{s-r-1} \binom{l_{n-1}-n+r}{r-1}
$$
Therefore
\begin{multline*}
\alpha(n-1;1,2,\ldots,i-2,i,\ldots,n-2,l_{n-2},l_{n-1}) \\
= \sum_{s=1}^{n} \sum_{t=1}^{n} \sum_{r=1}^{s-1} (-1)^{n+r+1}
|\overline{A}|_{n,s,i,t} \binom{-2n+s+1}{s-r-1} \binom{l_{n-1}-n+r}{r-1} \binom{l_{n-2}-n+t-1}{t-2}.
\end{multline*}
Consequently,
\begin{multline*}
\left.  \delta_{l_1} \delta_{l_2} \cdots \delta_{l_{n-3}}
\alpha_1(n-1;l_1,\ldots,l_{n-1})
\right|_{(l_1,\ldots,l_{n-3})=(1,2,\ldots,i-2,i,\ldots,n-2)} \\
= \sum_{s=1}^{n} \sum_{t=1}^{n} \sum_{r=1}^{s-1} (-1)^{r+1+n}
|\overline{A}|_{n,s,i,t} \binom{-2n+s+1}{s-r-1} \binom{l_{n-1}-n+r+1}{r} \binom{l_{n-2}-n+t}{t-1}.
\end{multline*}
By the Chu--Vandermonde summation
$$
\sum_{r=1}^{s-1} \binom{-2n+s+1}{s-r-1} \binom{l_{n-1}-n+r+1}{r} (-1)^{r+1} =
(-1)^{s} \left( \binom{n+l_{n-1}-1}{s-1} - \binom{2n-3}{s-1} \right),
$$
this simplifies to
\begin{multline*}
\left.  \delta_{l_1} \delta_{l_2} \cdots \delta_{l_{n-3}}
\alpha_1(n-1;l_1,\ldots,l_{n-1})
\right|_{(l_1,\ldots,l_{n-3})=(1,2,\ldots,i-2,i,\ldots,n-2)} \\
= \sum_{s=1}^{n} \sum_{t=1}^{n} (-1)^{n+s}
|\overline{A}|_{n,s,i,t} \left( \binom{n+l_{n-1}-1}{s-1} - \binom{2n-3}{s-1} \right) \binom{l_{n-2}-n+t}{t-1}.
\end{multline*}

\medskip

This implies
\begin{multline*}
\sum_{(l_{n-2},l_{n-1})}^{(n-1,k_{n-1},k_n)} \alpha(n-1;1,2,\ldots,i-2,i,\ldots,
n-2,l_{n-2},l_{n-1}) \\
= \sum_{(l_{n-2},l_{n-1})}^{(n-1,k_{n-1},k_n)} \left. \delta_{l_1} \delta_{l_2} \cdots
\delta_{l_{n-1}} \alpha_1(n-1;l_1,\ldots,l_{n-1}) \right|_{(l_1,\ldots,l_{n-3})=
(1,\ldots,i-2,i,\ldots,n-2)} \\
=  \left. \delta_{l_1} \delta_{l_2} \cdots
\delta_{l_{n-3}} \alpha_1(n-1;l_1,\ldots,l_{n-3},k_{n-1},k_n) \right|_{(l_1,\ldots,l_{n-3})=
(1,\ldots,i-2,i,\ldots,n-2)} \\
- \left. \delta_{l_1} \delta_{l_2} \cdots
\delta_{l_{n-3}} \alpha_1(n-1;l_1,\ldots,l_{n-3},n-2,k_n) \right|_{(l_1,\ldots,l_{n-3})=
(1,\ldots,i-2,i,\ldots,n-2)} \\
+ \left. \delta_{l_1} \delta_{l_2} \cdots
\delta_{l_{n-3}} \alpha_1(n-1;l_1,\ldots,l_{n-3},n-2,k_{n-1}-1) \right|_{(l_1,\ldots,l_{n-3})=
(1,\ldots,i-2,i,\ldots,n-2)} \\
= \sum_{s=1}^{n} \sum_{t=1}^{n} (-1)^{n+s} |\overline{A}|_{n,s,i,t}
\left( \binom{n+k_{n}-1}{s-1} - \binom{2n-3}{s-1} \right) \binom{k_{n-1}-n+t}{t-1} \\
+ [i=n] \sum_{s=1}^{n} (-1)^{n+s} A_{n-1,s} \left( \binom{n+k_{n-1}-2}{s-1} - \binom{n+k_n-1}{s-1} \right),
\end{multline*}
as $|\overline{A}|_{n,s,i,1} = [i=n] A_{n-1,s}$.

\medskip

\medskip

Therefore and by \eqref{shift}, we have
\begin{multline*}
\alpha(n;1,2,\ldots,n-3,k_{n-2},k_{n-1},k_n)
= \sum_{i=2}^{n-1} \sum_{s=1}^{n} \sum_{t=1}^{n} (-1)^{n+s} |\overline{A}|_{n,s,i,t} \\
\times \left( \binom{k_{n-2}-n+i+2}{i} \binom{k_{n-1} - n +
t}{t-1}
 \left(  \binom{n+k_n-1}{s-1} -  \binom{2n-3}{s-1} \right) \right. \\
- \binom{k_{n-2}-n+i+1}{i} \binom{k_{n-1} - n + t + 1}{t-1}  \left(
 \binom{n+k_n-1}{s-1} -
 \binom{2n-3}{s-1} \right) \\
- \binom{k_{n-2}-n+i+1}{i} \binom{k_{n-1} - n + t}{t-1} \left(
 \binom{n+k_n}{s-1} -
 \binom{2n-3}{s-1} \right) \\
\left. + \binom{k_{n-2}-n+i}{i} \binom{k_{n-1} - n + t + 1}{t-1} \left(
 \binom{n+k_n}{s-1} -  \binom{2n-3}{s-1} \right) \right) + R.
\end{multline*}

\subsection{The numbers $\overline{\overline{\underline{A_{n,s,i,t}}}}$ come into the play} On the other hand, by \eqref{cycle} and \eqref{expansion},
\begin{multline}
\label{2topbottom}
\alpha(n;1,2,\ldots,n-3,k_{n-2},k_{n-1},k_n) 
=
(-1)^{n-1} \alpha(n;k_n+n,1,2,\ldots,n-2,k_{n-2},k_{n-1}) \\
= (-1)^{n-1} \alpha(n;k_n+n+1,2,3,\ldots, n-2,k_{n-2}+1,k_{n-1}+1) \\
\sum_{s=1}^{n} \sum_{i=1}^{n} \sum_{t=1}^{n} \overline{\overline{\underline{A_{n,s,i,t}}}} (-1)^{n+s}
\binom{k_{n-2}-n+1+i}{i-1} \binom{k_{n-1}-n+1+t}{t-1} \binom{k_n+n-1}{s-1} 
\end{multline}

\medskip

It is not hard to
see that
{ \begin{multline}
\label{basis}
\left( \binom{k_{n-2}-n+i+2}{i} \binom{k_{n-1} - n + t}{t-1}
 \left(  \binom{n+k_n-1}{s-1} -  \binom{2n-3}{s-1} \right) \right. \\
- \binom{k_{n-2}-n+i+1}{i} \binom{k_{n-1} - n + t + 1}{t-1}  \left(
 \binom{n+k_n-1}{s-1} -
 \binom{2n-3}{s-1} \right) \\
- \binom{k_{n-2}-n+i+1}{i} \binom{k_{n-1} - n + t}{t-1} \left(
 \binom{n+k_n}{s-1} -
 \binom{2n-3}{s-1} \right) \\
\left. + \binom{k_{n-2}-n+i}{i} \binom{k_{n-1} - n + t + 1}{t-1} \left(
 \binom{n+k_n}{s-1} -  \binom{2n-3}{s-1} \right) \right) 
\end{multline}}
is equal to
\begin{multline*}
\left(  \binom{k_{n-2} - n +2 + i}{i} \binom{k_{n-1}-n+t}{t-2} - \binom{k_{n-2} - n + 1 + i}{i-1} \binom{k_{n-1}-n+1+t}{t-1} \right. \\
\left. + \binom{k_{n-2} - n + i}{i-2} \binom{k_{n-1}-n+1+t}{t-1} - \binom{k_{n-2} - n + 1 + i}{i-1} \binom{k_{n-1}-n+t}{t-2} \right) \\
\times \binom{n+k_n-1}{s-2} \\
+  \left(
\binom{k_{n-2} - n + i}{i-2} \binom{k_{n-1}-n+1+t}{t-1} - \binom{k_{n-2} - n + 1 + i}{i-1} \binom{k_{n-1}-n+t}{t-2} \right) \\
\times \binom{n+k_n-1}{s-1} \\
+ \left( \binom{k_{n-2} - n + 1 + i}{i-1} \binom{k_{n-1}-n+t}{t-2} - \binom{k_{n-2} - n + i}{i-2} \binom{k_{n-1}-n+1+t}{t-1}\right) \\
\times \binom{2n-3}{s-1}.
\end{multline*}
For the sake of brevity we set $f_{m,i}= \binom{k_m-n+1+i}{i-1}$ if $m=n-2,n-1$ and 
$f_{n,i} = \binom{n+k_n}{i-1}$.
Using this, we have 
\begin{multline*}
\alpha(n;1,2,\ldots,n-3,k_{n-2},k_{n-1},k_n) = \sum_{i=2}^{n-1} \sum_{s=1}^{n} \sum_{t=1}^{n}
|\overline{A}|_{n,s,i,t} (-1)^{n+s} \\
\times \left( f_{n-2,i+1} f_{n-1,t-1} f_{n,s-1} - f_{n-2,i} f_{n-1,t} f_{n,s-1} +
f_{n-2,i-1} f_{n-1,t} f_{n,s-1} - f_{n-2,i} f_{n-1,t-1} f_{n,s-1} \right. \\
\left. + f_{n-2,i-1} f_{n-1,t} f_{n,s} - f_{n-2,i} f_{n-1,t-1} f_{n,s} +
f_{n-2,i} f_{n-1,t-1} f_{n,1} \binom{2n-3}{s-1} - f_{n-2,i-1} f_{n-1,t} f_{n,1} \binom{2n-3}{s-1} \right) \\
+ \sum_{t=1}^{n} \left( (A_{n-1,t} - A_{n-1,t-1}) f_{n-2,n-2} f_{n-1,t} f_{n,1} 
+ (A_{n-1,t} - A_{n-1,t+1}) f_{n-2,n-1} f_{n-1,t} f_{n,1} \right) \\
+ \sum_{s=1}^{n} \left( (-1)^{n+s} (A_{n-1,s-1} - A_{n-1,s}) f_{n-2,n-2} f_{n-1,1} f_{n,s} 
+ (-1)^{n+s} A_{n-1,s} f_{n-2,n-1} f_{n-1,1} f_{n,s} \right).
\end{multline*}
After shifting certain indices, rearranging sums and using the fact that
$|\overline{A}|_{n,s,i,t}=0$ if either $s=n+1$ or $t=n+1$ or $i=0$ or ($i=1$ and $s > 1$) or  ($i=n$ and $t > 1$), we 
see that this is equal to
\begin{multline*}
\sum_{i=1}^{n} \sum_{s=1}^{n} \sum_{t=1}^{n} (-1)^{n+s} f_{n-2,i} f_{n-1,t} f_{n,s} \left( -|\overline{A}|_{n,s+1,i-1,t+1} + |\overline{A}|_{n,s+1,i,t} - |\overline{A}|_{n,s+1,i+1,t} + |\overline{A}|_{n,s+1,i,t+1} \right. \\ \left. + |\overline{A}|_{n,s,i+1,t} - |\overline{A}|_{n,s,i,t+1} \right) \\
+ \sum_{s=1}^{n} A_{n-1,s+1} (-1)^{n+s+1} f_{n-2,n} f_{n-1,n} f_{n,s} 
+ \sum_{s=1}^{n} A_{n-1,s+1} (-1)^{n+s} f_{n-2,n-1} f_{n-1,1} f_{n,s} \\
+ \sum_{s=1}^{n}  \left( A_{n-1,s-1} - A_{n-1,s} \right) (-1)^{n+s} f_{n-2,n-2} f_{n-1,1} f_{n,s} 
+ \sum_{t=1}^{n} \left( A_{n-1,t} - A_{n-1,t-1} \right) f_{n-2,n-2} f_{n-1,t} f_{n,1} \\
+ \sum_{t=1}^{n} \left( A_{n-1,t} - A_{n-1,t+1} \right) f_{n-2,n-1} f_{n-1,t} f_{n,1} + \sum_{i=1}^{n} \sum_{t=1}^{n} \sum_{j=1}^{n} |\overline{A}|_{n,j,i,t+1} (-1)^{n+j} \binom{2n-3}{j-1} f_{n-2,i} f_{n-1,t} f_{n,1} \\
+ \sum_{i=1}^{n} \sum_{t=1}^{n} \sum_{j=1}^{n} |\overline{A}|_{n,j,i+1,t} (-1)^{n+j+1} \binom{2n-3}{j-1} f_{n-2,i} f_{n-1,t} f_{n,1} \\
+ \sum_{j=1}^{n} A_{n-1,j} (-1)^{n+j} \binom{2n-3}{j-1} f_{n-2,n-1} f_{n-1,1} f_{n,1}.
\end{multline*}
We compare coefficients and conclude that 
\begin{multline*}
\overline{\overline{\underline{A_{n,s,i,t}}}} = -|\overline{A}|_{n,s+1,i-1,t+1} + |\overline{A}|_{n,s+1,i,t} - |\overline{A}|_{n,s+1,i+1,t} 
+ |\overline{A}|_{n,s+1,i,t+1} + |\overline{A}|_{n,s,i+1,t} - |\overline{A}|_{n,s,i,t+1} \\
-  [i=n][t=1] A_{n-1,s+1}  + [i=n-1][t=1] A_{n-1,s+1} 
+ [i=n-2] [t=1] \left( A_{n-1,s-1} - A_{n-1,s} \right) \\
+  [i=n-2] [s=1] (-1)^{n} \left( A_{n-1,t-1} - A_{n-1,t} \right)  
+ [i=n-1] [s=1] (-1)^{n} \left( A_{n-1,t+1} - A_{n-1,t} \right)  \\
+ [s=1] \sum_{j=1}^{n} (-1)^{j+1} |\overline{A}|_{n,j,i,t+1} \binom{2n-3}{j-1} 
+ [s=1] \sum_{j=1}^{n} (-1)^{j} |\overline{A}|_{n,j,i+1,t} \binom{2n-3}{j-1} \\
+ [i=n-1] [t=1] [s=1] \sum_{j=1}^{n} (-1)^{j+1}  A_{n-1,j} (-1)^{j+1} \binom{2n-3}{j-1}
\end{multline*}
The last expression can be simplified.
Observe that 
$$
\binom{2n-3}{j-1} = \binom{2n-4}{j-1} + \binom{2n-4}{j-2} = 
\binom{2n-4}{j-1} + \binom{2n-5}{j-2} + \binom{2n-5}{j-3} = \ldots =
\sum_{s=1}^{j} \binom{2n-3-s}{j-s}.
$$
Therefore
\begin{multline*}
\sum_{j=1}^{n} A_{n-1,j} (-1)^{j+1} \binom{2n-3}{j-1}
= \sum_{j=1}^{n} \sum_{s=1}^{j} A_{n-1,j} (-1)^{j+1} \binom{2n-3-s}{j-s}  \\
\sum_{s=1}^{n} \sum_{j=s}^{n} A_{n-1,j} (-1)^{j+1} \binom{2n-3-s}{j-s} = \sum_{s=1}^{n} (-1)^{n} A_{n-1,s} = (-1)^{n} A_{n-1}.
\end{multline*}
For $s > 1$ and $t > 1$, this gives
$$
\overline{\overline{\underline{A_{n,s,i,t}}}} = |\overline{A}|_{n,s+1,i,t} + |\overline{A}|_{n,s+1,i,t+1} + |\overline{A}|_{n,s,i+1,t} - |\overline{A}|_{n,s+1,i+1,t} - |\overline{A}|_{n,s,i,t+1} - |\overline{A}|_{n,s+1,i-1,t+1}
$$
and concludes the proof of Theorem~\ref{main}.

\medskip

\subsection{A final remark}

A question that remains is whether or not the linear relations \eqref{stroganov}, \eqref{doublyilse} and \eqref{triplytheo} have more refined generalizations. 
The fact that the coefficients of the various refined enumeration numbers in these relations are  always $1$ or $-1$ makes the relations more attractive and presumably more accessible for bijective proofs.

\medskip

For instance, if we reconsider \eqref{doublyilse} then we see that for all 
$i,t$ with $1 \le i < t \le n$ we have 
\begin{equation}
\label{doublyilseversion}
\overline{A}|_{n,i,t} = \overline{A}|_{n,i-1,t+1} + A||_{n,i,t} 
\end{equation}
where $A||_{n,i,t}$ is a version of $\overline{\overline{A}}_{n,i,t}$ in which the respective objects have gone through a rotation of $90$ degrees. Let 
$\overline{{\mathcal A}}|_{n,i,t}$ and ${\mathcal A}||_{n,i,t}$ denote  the set of objects that are counted by $\overline{A}|_{n,i,t}$ and $A||_{n,i,t}$. The relation \eqref{doublyilseversion} clearly suggests (at least to a ``bijective combinatorialist'') that there is a natural
decomposition of $\overline{{\mathcal A}}|_{n,i,t}$ into two sets such that 
one is in bijective relation with $\overline{\mathcal A}|_{n,i-1,t+1}$ and the other with ${\mathcal A}||_{n,i,t}$.

\medskip

Now, \eqref{triplytheo} could certainly assist in finding this decomposition and the bijections. In particular, it gives information on the role of the position of the unique $1$ in the leftmost column of the objects that are involved. Indeed, if we rewrite \eqref{triplytheo} as follows
\begin{equation}
\label{triplytheoversion}
|\overline{A}|_{n,s+1,i,t}
+ |\overline{A}|_{n,s+1,i,t+1} + |\overline{A}|_{n,s,i+1,t}
= |\overline{A}|_{n,s+1,i-1,t+1} + |A||_{n,s,i,t} + |\overline{A}|_{n,s,i,t+1} 
+ |\overline{A}|_{n,s+1,i+1,t}
\end{equation}
then we first of all notice that \eqref{doublyilseversion} follows after summing over all $s$ and cancelling $\overline{A}|_{n,i,t+1}$ and $\overline{A}|_{n,i+1,t}$. Suppose we have a bijective explanation of \eqref{triplytheoversion}: let $X^s_1 = |\overline{\mathcal A}|_{n,s+1,i,t}, 
X^s_2 = |\overline{\mathcal A}|_{n,s+1,i,t+1}, 
X^s_3 = |\overline{\mathcal A}|_{n,s,i+1,t}$ and 
$Y^s_1 = |\overline{\mathcal A}|_{n,s+1,i-1,t+1}, 
Y^s_2 = |{\mathcal A}||_{n,s,i,t}, Y^s_3 = |\overline{\mathcal A}|_{n,s,i,t+1}, 
Y^s_4 = |\overline{\mathcal A}|_{n,s+1,i+1,t}$. For all $i \in \{1,2,3\}$, let 
$X^s_{i,1} \cup X^s_{i,2} \cup X^s_{i,3} \cup X^s_{i,4} = X^s_i$ be a decomposition of $X^s_i$ and, for all $j \in \{1,2,3,4\}$, let $Y^s_{j,1} \cup Y^s_{j,2} \cup Y^s_{j,3} = Y^s_j$ be a decomposition of $Y^s_j$ and, for $(i,j) \in \{1,2,3\} \times \{1,2,3,4\}$, 
bijections $\phi^s_{i,j} : X^s_{i,j} \to Y^s_{j,i}$. (Note that although this appears to be more complicated compared to finding a bijective explanation for \eqref{doublyilseversion}, it is likely that this is in fact easier as the sets are smaller.)

\medskip

Given these decompositions and bijections, in order to construct a bijective explanation for 
\eqref{doublyilseversion} one may proceed as follows: let $A \in \overline{\mathcal A}|_{n,i,t}$. Then there is an $s$ such that $A \in X^s_1$. If $A \in X^s_{1,1}$ then $\phi_{1,1}$ maps 
$A$ to an element in  $\overline{\mathcal A}|_{n,i-1,t+1}$ and if 
$A \in X^s_{1,2}$ then $\phi_{1,2}$ maps $A$ to an element in 
${\mathcal A}||_{n,i,t}$  and we are done. Otherwise $\phi_{1,3}$ or $\phi_{1,4}$ maps 
$A$ to $Y^s_3$ or $Y^s_4$ respectively. However, via $Y^s_3=X^{s-1}_2$ and 
$Y^s_4 = X^{s+1}_3$ it is possible to change to the ``$X$-side'' of the equation again. Now, if the image lies in $X^{s-1}_{2,1} \cup X^{s-1}_{2,2}$ or 
$X^{s+1}_{3,1} \cup X^{s+1}_{3,2}$ then either $\phi^{s-1}_{2,1} \cup \phi^{s-1}_{2,2}$
or $\phi^{s+1}_{3,1} \cup \phi^{s+1}_{3,2}$ map this image to an element in 
$\overline{\mathcal A}|_{n,i-1,t+1} \cup {\mathcal A}||_{n,i,t}$. If not, one has to iterate this process.  As the sets are finite, this process has to 
terminate with an element in $\overline{\mathcal A}|_{n,i-1,t+1} \cup {\mathcal A}||_{n,i,t}$. All told, we are confronted with the rather odd fact that it may be more natural to find a 
bijective proof for
$$
\overline{A}|_{n,i,t} + \overline{A}|_{n,i,t+1} + \overline{A}|_{n,i+1,t} = \overline{A}|_{n,i-1,t+1} + A||_{n,i,t} + \overline{A}|_{n,i,t+1} + \overline{A}|_{n,i+1,t}
$$
than for \eqref{doublyilseversion}.

\bigskip \noindent
\textsc{
\!\!Ilse Fischer\\
Fakult\"at f\"ur Mathematik, Universit\"at Wien \\
1090 Wien, Austria \\
}
\texttt{Ilse.Fischer@univie.ac.at}

\end{document}